\documentclass{article}
\usepackage{amssymb,amsmath,amsthm,titlesec}
\usepackage[all]{xy}

\titleformat{\section}[block]
  {\sc\large\filcenter}
  {\S\thesection.}{.5em}{}

\theoremstyle{plain}
\newtheorem{theorem}{Theorem}[section]

\theoremstyle{plain}
\newtheorem{corollary}[theorem]{Corollary}

\theoremstyle{plain}
\newtheorem{prop}[theorem]{Proposition}

\theoremstyle{plain}
\newtheorem{lemma}[theorem]{Lemma}

\theoremstyle{definition}
\newtheorem{defn}[theorem]{Definition}

\theoremstyle{definition}
\newtheorem{example}[theorem]{Example}

\theoremstyle{definition}
\newtheorem{subsec}[theorem]{\S\hspace{-0.02in}}

\def\t{\textrm}
\def\ds{\displaystyle}

\def\bb{\mathbb}
\def\mb{\mathbf}
\def\mc{\mathcal}
\def\mf{\mathfrak}
\def\beq{\begin{eqnarray}}
\def\eeq{\end{eqnarray}}
\def\beqq{\begin{eqnarray*}}
\def\eeqq{\end{eqnarray*}}

\def\d{\partial}
\def\dbar{\bar{\d}}
\def\T{\mc{T}}
\def\e{\epsilon}
\def\n{\mathbf n}
\def\m{\mathbf m}
\def\bs{\boldsymbol}
\def\wt{\widetilde}
\def\vs{\vspace{0.01in}}
\def\vss{\vspace{0.02in}}
\def\ph{\phantom}
\def\bvarphi{{\bs{\varphi}}}
\def\ve{\varepsilon}
\def\CDO{\mc{D}^{\mathrm{ch}}}
\def\Tr{\mathrm{Tr}\,}
\def\Str{\mathrm{Str}\,}
\def\tnabla{\nabla^t}

\title{\sf Chiral Differential Operators on Supermanifolds}
\author{\sf Pokman Cheung}
\date{\sf June 22, 2011}

\begin{document}

\maketitle

\begin{abstract}
The first part of this paper provides a new description of
chiral differential operators (CDOs) in terms of global
geometric quantities.
The main result is a recipe to define all sheaves of CDOs on
a smooth cs-manifold;
its ingredients consist of an affine connection $\nabla$ and
an even $3$-form that trivializes $p_1(\nabla)$.
With $\nabla$ fixed, two suitable $3$-forms define isomorphic
sheaves of CDOs if and only if their difference is exact.
Moreover, conformal structures are in one-to-one correspondence
with even $1$-forms that trivialize $c_1(\nabla)$.

Applying our work in the first part, we construct what may be
called ``chiral Dolbeault complexes'' of a complex manifold $M$,
and analyze conditions under which these differential vertex
superalgebras admit compatible conformal structures or extra
gradings (fermion numbers).
When $M$ is compact, their cohomology computes (in various
cases) the Witten genus, the two-variable elliptic genus and
a spin$^c$ version of the Witten genus.
This part contains some new results as well as provides
a geometric formulation of certain known facts from
the study of holomorphic CDOs and $\sigma$-models.
\end{abstract}

\section{Introduction}

In physics, the study of a type of quantum field theory called
\emph{$\sigma$-models} has inspired many important insights in topology
and geometry.
The theory of elliptic genera is an example.
In particular, associated to any compact, string manifold
\footnote{
Let $\lambda\in H^4(B\t{Spin};\bb{Z})\cong\bb{Z}$ be the generator
such that $2\lambda=p_1$.
This defines a characteristic class $\lambda(\cdot)$ for spin vector
bundles.
A spin manifold $M$ is said to be \emph{string} if $\lambda(TM)=0$.
Moreover, a \emph{string structure} on $M$ is a ``trivialization of
$\lambda(TM)$'', i.e.~a homotopy class of liftings of the classifying map
$M\rightarrow B\t{Spin}$ along the homotopy fiber of
$\lambda:B\t{Spin}\rightarrow K(\bb{Z},4)$.}
$M$ is a $\sigma$-model whose ``partition function'' equals, up to
a constant factor, the formal power series
\beqq
W(M)=\int_M\hat{A}(TM)\,
  ch\bigg(\bigotimes_{n=1}^\infty\t{Sym}_{q^n}(TM\otimes\bb{C})\bigg)
  \cdot\prod_{n=1}^\infty(1-q^n)^{\dim M}
\eeqq
known as the \emph{Witten genus} of $M$.~\cite{Witten.ell,Witten.index}
Similarly, associated to any compact, spin manifold $M$ is another
$\sigma$-model, which gives rise to the formal power series
\beqq
\t{Och}(M)=\int_M L(TM)\,
  ch\bigg(\bigotimes_{n=1}^\infty\t{Sym}_{q^n}(TM\otimes\bb{C})\otimes
  \bigotimes_{n=1}^\infty\wedge_{q^n}(TM\otimes\bb{C})\bigg)
  \cdot\prod_{n=1}^\infty\left(\frac{1-q^n}{1+q^n}\right)^{\dim M}
\eeqq
known as the \emph{Ochanine elliptic genus} of
$M$.~\cite{Ochanine,Witten.ell}
The physical interpretation of these topological invariants have
led to predictions that are not immediately clear from 
the mathematical point of view.
Even though many of them have since been verified,
e.g.~\cite{Zagier,Bott.Taubes},
a complete, geometric understanding of elliptic genera has
yet to emerge.
The latter probably requires to some extent a mathematical
framework for $\sigma$-models.

Sheaves of vertex algebras provide a mathematical approach to
$\sigma$-models.
Important constructions along this line include
the \emph{chiral de Rham complex} and, more generally,
sheaves of \emph{chiral differential operators}, or
\emph{CDOs}.~\cite{MSV,GMS1}
In particular, a complex manifold $M$ admits a sheaf of
holomorphic CDOs $\CDO_M$ with a conformal structure if and
only if $c_1^{\t{hol}}(TM)=c_2^{\t{hol}}(TM)=0$;
\footnote{
See Definition \ref{hol.Chern}.}
notice that $M$ as a spin$^c$ manifold admits a string structure
if and only if $c_1(TM)=c_2(TM)=0$.
Furthermore, if $M$ is compact
\beqq
\t{char}\,H^*(M,\CDO_M)=W(M)\cdot(\t{a constant factor})
\eeqq
suggesting a connection between $\CDO_M$ and the $\sigma$-model
underlying the Witten genus.
In fact, physicists have recognized a connection between CDOs
and $\sigma$-models of various
flavors.~\cite{Kapustin,Witten.CDO,Tan}
More recently, a new construction of the Witten genus has
been given under a systematic mathematical framework for
perturbative quantum field theory.~\cite{Costello}

The first goal of this paper is to provide a new description
of CDOs using global geometric quantities and the language
of cs-manifolds, i.e.~supermanifolds equipped with
$\bb{C}$-valued functions.
The algebra of \emph{smooth} CDOs on $\bb{R}^{p|q}$ is
the smooth analogue of the conformal vertex superalgebra
$(\beta\gamma)^{\otimes p}\otimes(bc)^{\otimes q}$
(\S\ref{sheafCDO.Rpq}, Proposition \ref{conformal});
its behavior under a change of coordinates, first computed in \cite{GMS1},
are restated here in more geometric terms
(\S\ref{sec.CDO.coorchange}, Proposition \ref{conformal.coorchange}).
The notions of a sheaf of CDOs and its conformal structures are then
generalized from $\bb{R}^{p|q}$ to a general cs-manifold $\mb{M}$ in
a natural way (Definition \ref{sheafCDO.defn}).
After dealing with some technical issues (Lemmas \ref{globalCDO.lemma1},
\ref{globalCDO.lemma2}), we prove the main result on the global construction
of CDOs (Theorem \ref{thm.globalCDO}).
Namely, given an affine connection $\nabla$ and an even $3$-form $H$ that
satisfies
\beqq
dH=\Str(R\wedge R)
\eeqq
where $R=\t{curv}(\nabla)$, there is a recipe to define
a sheaf of CDOs $\CDO_{\mb{M},\nabla,H}$, and this recipe
yields essentially all sheaves of CDOs on $\mb{M}$.
Moreover, conformal structures on $\CDO_{\mb{M},\nabla,H}$
are in one-to-one correspondence with even $1$-forms
$\omega$ that satisfy
\beqq
d\omega=\Str R.
\eeqq
To classify these objects, we also prove that, with $\nabla$ fixed,
two suitable $3$-forms $H,H'$ define isomorphic sheaves of CDOs if and only
if $H-H'$ is exact (Theorem \ref{thm.globalCDO.iso}).
In contrast to \cite{GMS1}, our description of CDOs does not
rely on a choice of coordinate charts or other local data.
For the special case of the chiral de Rham complex, in
which both $H$ and $\omega$ are trivial (Example \ref{cdR}),
an invariant description has also been given in \cite{BHS}.
The formulation of CDOs developed here has been applied
e.g.~to study how to lift a Lie group action on a manifold to
a ``formal loop group action'' on CDOs.~\cite{Cheung}

In the rest of the paper, we apply our work in the first part
to construct what may be called ``chiral Dolbeault complexes.''
Let $M$ be a complex manifold and $E\rightarrow M$
a holomorphic vector bundle.
The Dolbeault complex of $M$ valued in $\wedge^*E^\vee$
is identified with the smooth functions on the cs-manifold
\beqq
\mb{M}=\Pi(\overline{TM}\oplus E)
\eeqq
under the action of an odd vector field $Q$ that
satisfies $Q^2=0$ (\S\ref{sec.Dolbcs}).
This motivates us to construct a sheaf of CDOs
$\CDO_{\mb{M},\nabla,H}$ on $\mb{M}$ (\S\ref{sec.DolbCDO}),
and study the condition under which the supersymmetry $Q$
lifts to one on CDOs, i.e.~an odd derivation
$\hat Q$ on $\CDO_{\mb{M},\nabla,H}$ that satisfies $\hat Q^2=0$
(Theorem \ref{thm.DolbQ}, Proposition \ref{prop.obstructions}).
At the same time we also analyze the condition for $\hat{Q}$ to
respect a conformal structure.
Moreover, if one or both of the line bundles $\det TM$, $\det E$
are flat, $\hat Q$ is compatible with certain gradings on
$\CDO_{\mb{M},\nabla,H}$ called fermion numbers
(\S\ref{sec.fermion}, Propositions \ref{prop.Q.Jr},
\ref{prop.Q.Jl}).
The sheaf of differential vertex superalgebras
\beqq
(\CDO_{\mb{M},\nabla,H},\hat{Q})
\eeqq
may be thought of as a Dolbeault resolution of holomorphic
CDOs on $\Pi E$, as well as a particular limit of
a $\sigma$-model.~\cite{Kapustin}
When $M$ is compact, its cohomology computes various elliptic
genera (Theorem \ref{thm.cDolb.chlgy}),
including the Witten genus in the case $E=0$ (Example \ref{E=0}),
a two-variable generalization of the Ochanine genus in the case
$E=TM$ (Example \ref{E=TM}), and a spin$^c$ version of the Witten
genus \cite{CHZ} in the cases $E=\det TM$ and
$E=(\det TM)^{\otimes 2}-\det TM$
(Examples \ref{E=l}, \ref{E=l2-l}).
Most of the results in this part are similar to and consistent
with what is known from the study of holomorphic CDOs and
$\sigma$-models, but our formulation may provide a new
geometric point of view.
On the other hand, the last two examples seem to be new.

The first appendix reviews the notion of vertex algebroids
(first introduced in \cite{GMS2}), their relation with
vertex algebras, and gives some examples.
Despite the rather complicated-looking definition, vertex
algebroids and their super analogues provide a convenient
tool in our study of CDOs.
In the second appendix, we construct affine connections on
cs-manifolds and obtain formulae that are needed in various
calculations with CDOs.

\vspace{0.05in}
\noindent
{\bf Conventions.}
For the definition of a vertex superalgebra, see \cite{Kac.VA,FB-Z}.
In this paper, every vertex superalgebra $V$ is graded
by non-negative integers called weights.
The notation $V_k$ means its component of weight $k$, and
$L_0$ denotes the weight operator, so that $L_0|_{V_k}=k$.

For the definition of a cs-manifold, see \cite{QFS.susy}.
Given a smooth cs-manifold $\mb{M}$, we always denote by
$C^\infty_{\mb{M}}$, $\T_{\mb{M}}$ and $\Omega^n_{\mb{M}}$
its sheaves of \emph{smooth} functions, vector fields and
$n$-forms;
when ``$\mb{M}$'' appears in parentheses instead of
the subscript, it means the corresponding spaces
of global sections.
Restricting $C^\infty_{\mb{M}}$ to an open subset
$U\subset\mb{M}^{\t{red}}$ defines a new cs-manifold, denoted
by $\mb{M}|_U$.
Square brackets are used for supercommutators between operators
of any parities, while ``$\Str$'' stands for the supertrace.
Notice that $\bb{R}^{p|q}$ is regarded as a cs-manifold in
this paper, namely
\beqq
C^\infty_{\bb{R}^{p|q}}
=C^\infty_{\bb{R}^p}\otimes\wedge^*(\bb{R}^q)\otimes\bb{C}.
\eeqq

\vspace{0.05in}
\noindent
{\bf Acknowledgements.}
This paper grew out of an effort to understand some of
the pioneering work by Vassily Gorbounov, Fyodor Malikov,
Vadim Schechtman and Arkady Vaintrob.
The author would like to thank Ralph Cohen, Haynes Miller,
Stephan Stolz and Peter Teichner for their continual interest
and encouragement.
He has also benefited from discussions with Vassily Gorbounov,
Fei Han, Qin Li and Meng-Chwan Tan.
He did most of the work in this paper under the generous
support of the Max-Planck-Institut f\"ur Mathematik, before
finishing it during a visit at the Hong Kong University of
Science and Technology.
He also gratefully acknowledges the current support of
an EPSRC grant.

\newpage

\setcounter{equation}{0}
\section{Chiral Differential Operators}

Sheaves of CDOs on a manifold were first studied in \cite{GMS1}.
This section provides an alternative construction of the smooth version 
using global geometric quantities.

\begin{subsec} \label{sheafCDO.Rpq}
{\bf The sheaf of CDOs on $\bb{R}^{p|q}$.}
Let $b^1,\ldots,b^p$ and $b^{p+1},\ldots,b^{p+q}$ be respectively 
the even and odd coordinates of $\bb{R}^{p|q}$.
The following notations are used
\beqq
\d_i=\frac{\d}{\d b^i}\,,\qquad
|\cdot|=\t{parity},\quad
\e_i=(-1)^{|b^i|},\quad
\e_{ij}=(-1)^{|b^i||b^j|}
\eeqq
and repeated indices are summed over (but not counting those from
$\e_i,\e_{ij}$).
Regard $\bb{R}^{p|q}$ as a smooth cs-manifold, namely
\beqq
\textstyle
C^\infty(\bb{R}^{p|q})=C^\infty(\bb{R}^p)\otimes
  \bigwedge(b^{p+1},\ldots,b^{p+q})\otimes\bb{C}.
\eeqq
Given an open set $W\subset\bb{R}^p$, let $\mb{W}=(\bb{R}^{p|q})|_W$.
Consider the vertex superalgebra $\CDO(\mb{W})$ constructed in 
\S\ref{smoothCDO}.
It is freely generated by a vertex superalgebroid 
\beqq
\big(C^\infty(\mb{W}),\Omega^1(\mb{W}),\T(\mb{W}),
  \ast^c,\{\;\}^c,\{\;\}^c_\Omega\big)
\eeqq
and, by the following result, equipped with a family of conformal elements
\beq \label{localCDO.conformal}
\nu^\omega:=\e_i\d_{i,-1}db^i+\frac{1}{2}\omega_{-2}\mb{1},\qquad
\omega\in\Omega^1(\mb{W}),\;|\omega|=\bar 0,\;d\omega=0.
\eeq
The assignment $W\mapsto\CDO_{p|q}(W):=\CDO(\mb{W})$ defines a sheaf of 
conformal vertex superalgebras on $\bb{R}^p$.
\end{subsec}

\begin{prop} \label{conformal}
The elements $\nu^\omega$ in (\ref{localCDO.conformal}) are conformal
in $\CDO(\mb{W})$ of central charge $2(p-q)$.
\end{prop}

\begin{proof}
First consider $\nu:=\e_i\d_{i,-1}db^i$.
Let us show that
\beqq
\t{(i)}\quad\nu_{(0)}=T,\;\nu_{(1)}=L_0\t{ on }C^\infty(\mb{W})\cup
  \{\d_k\}_{k=1}^{p+q}\qquad\qquad
\t{(ii)}\quad\nu_{(3)}\nu=p-q
\eeqq
The operators $\d_{i,n}$ for $i=1,\ldots,p+q$ and $n\in\bb{Z}$ commute 
with each other, because
\beqq
[\d_{i,n},\d_{j,m}]
=[\d_i,\d_j]_{n+m}
+\{\d_i,\d_j\}^c_{\Omega,n+m}
+n\{\d_i,\d_j\}^c_{n+m}=0.
\eeqq
Keeping this in mind, we compute the following for $f\in C^\infty(\mb{W})$ 
and $k=1,\ldots,p+q$
\begin{align*}
\nu_{(1)}f\; 
& =(db^i)_0\d_{i,0}f=0 \\
\nu_{(0)}f\;
& =\e_i\d_{i,-1}(db^i)_0 f+(db^i)_{-1}\d_{i,0}f
  =0+db^i\cdot\d_i f=df=Tf \\
\nu_{(2)}\d_k
& =(db^i)_1\d_{i,0}\d_k+(db^i)_0\d_{i,1}\d_k=0 \\
\nu_{(1)}\d_k
& =\e_i\d_{i,-1}(db^i)_1\d_k+(db^i)_0\d_{i,0}\d_k
  +(db^i)_{-1}\d_{i,1}\d_k
  =\e_i\d_{i,-1}db^i(\d_k)+0+0=\d_k \\
\nu_{(0)}\d_k
& =\e_i\d_{i,-2}(db^i)_1\d_k+\e_i\d_{i,-1}(db^i)_0\d_k
  +(db^i)_{-1}\d_{i,0}\d_k+(db^i)_{-2}\d_{i,1}\d_k \\
& =\e_i\d_{i,-2}db^i(\d_k)+0+0+0
  =\d_{k,-2}\mb{1}=T\d_k \\
\nu_{(3)}\nu\;
& =\e_i[\nu_{(3)},\d_{i,-1}]db^i
  =\e_i(\nu_{(0)}\d_i)_{(2)}db^i
  +3\e_i(\nu_{(1)}\d_i)_{(1)}db^i
  +3\e_i(\nu_{(2)}\d_i)_{(0)}db^i \\
& =-2\e_i\d_{i,1}db^i+3\e_i\d_{i,1}db^i+0
  =db^i(\d_i)=p-q
\end{align*}
This proves (i) and (ii).
Now notice that $[\nu_{(1)},f_0]=(\nu_{(0)}f)_{(0)}+(\nu_{(1)}f)_{(-1)}=0$, 
and also that both $\nu_{(0)},T$ are vertex superalgebra derivations 
commuting with $T$.
Then compute for $\alpha\in\Omega^1(\mb{W})$, $X\in\T(\mb{W})$
\begin{align*}
\nu_{(1)}\alpha\;
& =\nu_{(1)}(\alpha_k db^k)=\nu_{(1)}\alpha_{k,0}Tb^k
  =[\nu_{(1)},\alpha_{k,0}]Tb^k+\alpha_{k,0}[\nu_{(1)},T]b^k
  +\alpha_{k,0}T\nu_{(1)}b^k \\
& =0+\alpha_{k,0}\nu_{(0)}b^k+0=\alpha_k db^k=\alpha \\
\nu_{(0)}\alpha\;
& =\nu_{(0)}(\alpha_k db^k)=\nu_{(0)}\alpha_{k,0}Tb^k
  =T\alpha_{k,0}Tb^k=T\alpha \\
\nu_{(1)}X
& =\nu_{(1)}(X^k\d_k)=\nu_{(1)}X^k_0\d_k-\nu_{(1)}(X^k\ast^c\d_k) 
  =[\nu_{(1)},X^k_0]\d_k+X^k_0\nu_{(1)}\d_k-X^k\ast^c\d_k \\
& = 0+X^k_0\d_k-X^k\ast^c\d_k = X^k\d_k=X \\
\nu_{(0)}X
& =\nu_{(0)}(X^k\d_k)=\nu_{(0)}X^k_0\d_k-\nu_{(0)}(X^k\ast^c\d_k) 
  =TX^k_0\d_k-T(X^k\ast^c\d_k)=TX
\end{align*}
Hence (i) implies that $\nu_{(0)}=T$, $\nu_{(1)}=L_0$ also hold on 
$\Omega^1(\mb{W})$ and $\T(\mb{W})$.
This yields the commutation relations:
\beqq
\begin{array}{lll}
\lbrack\nu_{(0)},f_n]=(1-n)f_{n-1}\ph{aa} &
[\nu_{(0)},\alpha_n]=-n\alpha_{n-1}\ph{aa} &
[\nu_{(0)},X_n]=-nX_{n-1} \vss \\
\lbrack\nu_{(1)},f_n]=-nf_n &
[\nu_{(1)},\alpha_n]=-n\alpha_n &
[\nu_{(1)},X_n]=-nX_n
\end{array}
\eeqq
Since we also have $\nu_{(0)}\mb{1}=0=\nu_{(1)}\mb{1}$, the operators 
$\nu_{(0)},\nu_{(1)}$ satisfy respectively the defining relations of 
$T$ and $L_0$, i.e.~$\nu_{(0)}=T$, $\nu_{(1)}=L_0$ on the entire vertex 
superalgebra $\CDO(\mb{W})$.
By Lemma 3.4.5 of \cite{FB-Z}, this together with (ii) proves 
the proposition for $\nu$.

Let $\delta$ denote an even element of $\CDO(\mb{W})_2$.
Replacing $\nu$ by $\nu+\delta$ in the above arguments shows that 
$\nu+\delta$ is also conformal of the same central charge if
\beqq
\t{(i)}'\quad\delta_{(0)}=\delta_{(1)}=0\t{ on }C^\infty(\mb{W})\cup
  \{\d_k\}_{k=1}^{p+q}\qquad\qquad
\t{(ii)}'\quad\nu_{(3)}\delta+\delta_{(3)}\nu+\delta_{(3)}\delta=0
\eeqq
Suppose $\delta=\omega_{-2}\mb{1}=T\omega$, where 
$\omega\in\Omega^1(\mb{W})$.
Then $\delta_{(n)}=-n\omega_{n-1}$.
For $f\in C^\infty(\mb{W})$, $k=1,\ldots,p+q$
\begin{align*}
&\delta_{(0)}f=0=\delta_{(0)}\d_k \\
&\delta_{(1)}f=-\omega_0 f=0 \\
&\delta_{(1)}\d_k
  =-\omega_0\d_k=[\d_{k,-1},\omega_0]\mb{1}
  =L_{\d_k}\omega-d\omega(\d_k)=\iota_{\d_k}d\omega
\end{align*}
Hence (i)$'$ is satisfied if $d\omega=0$.
On the other hand, 
$[\nu_{(2)},f_0]=(\nu_{(0)}f)_{(1)}+2(\nu_{(1)}f)_{(0)}=-f_1$ implies
\begin{align*}
\nu_{(2)}\omega
  =\nu_{(2)}(\omega_k db^k)
& =\nu_{(2)}\omega_{k,0}Tb^k
  =-\omega_{k,1}Tb^k+\omega_{k,0}\nu_{(2)}Tb^k
  =0+2\omega_{k,0}\nu_{(1)}b^k=0 \\
\Rightarrow\qquad
\nu_{(3)}\delta+\delta_{(3)}\nu+\delta_{(3)}\delta
& =\nu_{(3)}\omega_{-2}\mb{1}-3\omega_2\nu-3\omega_2\omega_{-2}\mb{1} \\
& =[\nu_{(3)},\omega_{-2}]\mb{1}+3[\nu_{(-1)},\omega_2]\mb{1}
  -3[\omega_2,\omega_{-2}]\mb{1} \\ 
& =4(\nu_{(0)}\omega)_{(1)}\mb{1}+6\nu_{(2)}\omega=0
\end{align*}
so that (ii)$'$ holds.
This completes the proof of the proposition.
\end{proof}

{\it Remark.}
In the above proof, full details are shown in order to demonstrate 
the type of arguments involved in similar calculations.
Subsequent proofs will be given more briefly.

\begin{subsec} \label{sec.CDO.coorchange}  
{\bf Coordinate transformations of CDOs on $\bb{R}^{p|q}$.}
Let $\mb{W}$, $\mb{W}'$, $\mb{W}''$ be restrictions of $\bb{R}^{p|q}$ 
(as a cs-manifold) to open sets in $\bb{R}^p$. 
Suppose $\bvarphi:\mb{W}\rightarrow\mb{W}'$ is a diffeomorphism of 
cs-manifolds.
Recall the notations in \S\ref{MC.WZ} and Theorem \ref{CDO.iso}.
Given an even $2$-form $\xi$ on $\mb{W}$ with $d\xi=WZ_\bvarphi$, there 
is a corresponding isomorphism of vertex superalgebras
\beqq
\bvarphi^*_\xi:\CDO(\mb{W}')\rightarrow\CDO(\mb{W}).
\eeqq
For $f\in C^\infty(\mb{W}')$, $\alpha\in\Omega^1(\mb{W}')$ and 
$X\in\T(\mb{W}')$, we have
\beq \label{CDO.coorchange}
\bvarphi^*_\xi(f)=\bvarphi^*f,\qquad
\bvarphi^*_\xi(\alpha)=\bvarphi^*\alpha,\qquad
\bvarphi^*_\xi(X)=\bvarphi^*X+\Delta_{\bvarphi,\xi}(X).
\eeq
All isomorphisms between $\CDO(\mb{W}')$ and $\CDO(\mb{W})$ are of this 
form.
According to the result below, $\bvarphi^*_\xi$ permutes the conformal
elements (\ref{localCDO.conformal}).
If $\bvarphi':\mb{W}'\rightarrow\mb{W}''$ is another diffeomorphism of 
cs-manifolds and $\xi'$ is an even $2$-form on $\mb{W}'$ with 
$d\xi'=WZ_{\bvarphi'}$, then
\beqq
\bvarphi^*_\xi\circ\bvarphi^{\prime*}_{\xi'}
=(\bvarphi'\bvarphi)^*_\eta,\qquad
\eta=\xi+\bvarphi^*\xi'+\sigma_{\bvarphi',\bvarphi}.
\eeqq
\end{subsec}

\begin{prop} \label{conformal.coorchange}
Consider the isomorphism 
$\bvarphi^*_\xi:\CDO(\mb{W}')\rightarrow\CDO(\mb{W})$ described above.
For even closed $1$-forms $\omega$ on $\mb{W}'$, we have
\beqq
\bvarphi^*_\xi(\nu^\omega)=\nu^{\bvarphi^*\omega-\Str\theta_\bvarphi}.
\eeqq
\end{prop}

{\it Remark.}
Notice that $d\theta_\bvarphi=-\theta_\bvarphi\wedge\theta_\bvarphi$ implies
$\Str\theta_\bvarphi$ is closed.
Also, as a consistency check, it follows from 
$\theta_{\bvarphi'\bvarphi}=\theta_\bvarphi
+g_\bvarphi^{-1}\cdot\bvarphi^*\theta_{\bvarphi'}\cdot g_\bvarphi$ that 
$\Str\theta_{\bvarphi'\bvarphi}=\Str\theta_\bvarphi
+\bvarphi^*\Str\theta_{\bvarphi'}$.

\begin{proof}[Proof of Proposition \ref{conformal.coorchange}]
It suffices to consider the case $\omega=0$.
To simplify notations, let us write $g=g_\bvarphi$, $h=g^{-1}$, 
$\theta=\theta_\bvarphi$ and $\Delta=\Delta_{\bvarphi,\xi}$.
By (\ref{CDO.coorchange}), we have
\beqq
\bvarphi^*_\xi(\nu)
=\e_i\big(\bvarphi^*\d_i+\Delta(\d_i)\big)_{-1}(\bvarphi^*db^i)
=\e_{ik}(h^k_{\ph{i}i}\d_k)_{-1}(g^i_{\ph{i}\ell}db^\ell)
  +\e_i\Delta(\d_i)_{-1}(g^i_{\ph{i}\ell}db^\ell)
\eeqq
The first term above is computed as follows:
\begin{align*}
&\e_{ik}(h^k_{\ph{i}i}\d_k)_{-1}(g^i_{\ph{i}\ell}db^\ell) \\
=\;&\e_{ik}\Big(
  h^k_{\ph{i}i,-2}\,\d_{k,1}+h^k_{\ph{i}i,-1}\,\d_{k,0}
  +h^k_{\ph{i}i,0}\,\d_{k,-1}-(h^k_{\ph{i}i}\ast\d_k)_{-1}\Big)
  g^i_{\ph{i}\ell,0}\,db^\ell \\
=\;&\e_k h^k_{\ph{i}i,-2}\,g^i_{\ph{i}k}
  +\e_k h^k_{\ph{i}i,-1}\,dg^i_{\ph{i}k}
  +\e_k\d_{k,-1}db^k
  -\e_k(\d_k h^k_{\ph{i}i})_{-1}\,g^i_{\ph{i}\ell,0}\,db^\ell
  +\e_k(\d_k h^k_{\ph{i}i})_{-1}\,g^i_{\ph{i}\ell,0}\,db^\ell \\
=\;&\nu
  +\frac{1}{2}\,\Str\big((dh)_{-2}\,g\big)
  +\Str\big((dh)_{-1}\,dg\big)
= \nu
  -\frac{1}{2}\,\Str(\theta_{-2}\mb{1})
  -\frac{1}{2}\,\Str(\theta_{-1}\theta)
\end{align*}
Then we compute the second term above:
\begin{align*}
&\e_i\Delta(\d_i)_{-1}(g^i_{\ph{i}\ell}db^\ell) \\
=\;&\left(
  -\e_k\e_r\e_{ik}\e_{ir}\e_{kr}\d_r h^k_{\ph{i}i}\cdot\theta^r_{\ph{i}k}
  -\frac{1}{2}\e_{ik}\Str(\theta\otimes\theta)
    (h^k_{\ph{i}i}\d_k\otimes\,\t{-}\,)
  -\frac{1}{2}\e_{ik}\,\xi(h^k_{\ph{i}i}\d_k,\,\t{-}\,)\right)_{-1}
  \hspace{-0.03in}g^i_{\ph{i}\ell,0}\,db^\ell \\
=\;&\e_r\theta^r_{\ph{i}k,-1}\,h^k_{\ph{i}i,0}\,dg^i_{\ph{i}r}
  -\frac{1}{2}\e_k\e_r\Str(\theta\otimes\theta)(\d_k\otimes\d_r)_0
    (db^r)_{-1}db^k
  -\frac{1}{2}\e_k\e_r\xi(\d_k,\d_r)_0(db^r)_{-1}db^k \\
=\;&\Str(\theta_{-1}\theta)-\frac{1}{2}\Str(\theta_{-1}\theta)
  =\frac{1}{2}\Str(\theta_{-1}\theta)
\end{align*}
where we have used the graded symmetry of 
$(db^r)_{-1}db^k=b^r_{-1}b^k_{-1}\mb{1}$.
This yields
\[
\bvarphi^*_\xi(\nu)
=\nu-\frac{1}{2}(\Str\theta)_{-2}\mb{1}
=\nu^{-\Str\theta}. \qedhere
\]
\end{proof}

{\it Preparation.}
Given topological spaces $X,X'$, a presheaf $\mc{S}$ on $X$ and a presheaf 
$\mc{S}'$ on $X'$ valued in some category, let 
$(\varphi,\Phi):(X,\mc{S})\rightarrow(X',\mc{S}')$ denote the data
consisting of a continuous map $\varphi:X\rightarrow X'$ and a morphism of
presheaves $\Phi:\mc{S}'\rightarrow\varphi_*\mc{S}$ on $X'$.
Composition reads
$(\varphi',\Phi')\circ(\varphi,\Phi)
=(\varphi'\varphi,\,\varphi'_*\Phi\circ\Phi')$.
Recall the sheaf of vertex superalgebras $\CDO_{p|q}$ described in 
\S\ref{sheafCDO.Rpq}.

\begin{defn} \label{sheafCDO.defn}
A \emph{sheaf of CDOs on a smooth $(p|q)$-dimensional cs-manifold 
$\mb{M}=(M,C^\infty_{\mb{M}})$} is a sheaf of vertex superalgebras $\mc{V}$ 
on $M$ with the following properties: \vspace{-0.05in}
\begin{itemize}
\item[$\cdot$]
The weight-zero component is $\mc{V}_0=C^\infty_{\mb{M}}$.
\vspace{-0.08in}
\item[$\cdot$]
Given $x\in M$, there exist open sets $U\subset M$, $W\subset\bb{R}^p$
with $x\in U$, and an isomorphism between $(U,\mc{V}|_U)$ and 
$(W,\CDO_{p|q}|_W)$ as topological spaces equipped with sheaves of 
vertex superalgebras.
\vspace{-0.08in}
\end{itemize}
A \emph{conformal structure} on $\mc{V}$ is an element 
$\nu\in\mc{V}(M)_2$ such that, under each isomorphism postulated
above, $\nu|_U\in\mc{V}(U)$ corresponds to one of the conformal
elements $\nu^\omega\in\CDO_{p|q}(W)$ described in
(\ref{localCDO.conformal}).
\end{defn}

{\it Remark.}
For example, $\CDO_{p|q}$ is a sheaf of CDOs on $\bb{R}^{p|q}$ with 
a family of conformal structures $\nu^\omega$.
While a general sheaf of CDOs is locally isomorphic to $\CDO_{p|q}$, 
the latter has up to this point been defined using coordinates in
a manifest way (see \S\ref{sheafCDO.Rpq} and appendix \S\ref{app.VAoid}).
The geometric data required to globalize the construction is the main
content of Theorem \ref{thm.globalCDO}. \vspace{0.05in}

{\it Preparation.}
The sheaves of smooth functions, $1$-forms and vector fields on a smooth
cs-manifold $\mb{M}$ form a sheaf of extended Lie superalgebroids 
$(C^\infty_{\mb{M}},\Omega^1_{\mb{M}},\T_{\mb{M}})$ using the usual 
differential on functions, Lie brackets on vector fields, Lie derivations
on functions and $1$-forms by vector fields, and pairing between $1$-forms
and vector fields.

\begin{lemma} \label{globalCDO.lemma1}
Let $(\varphi,\Phi):(U,\mc{V}|_U)\rightarrow(W,\CDO_{p|q}|_W)$ be 
an isomorphism as postulated in Definition \ref{sheafCDO.defn}.
Also let $\mb{U}=\mb{M}|_U$, $\mb{W}=(\bb{R}^{p|q})|_W$.

(a) The data determine a diffeomorphism of cs-manifolds 
$\bvarphi:\mb{U}\rightarrow\mb{W}$.
The presheaf (in fact, sheaf) of extended Lie superalgebroids associated
to $\mc{V}|_U$ can be identified with 
$(C^\infty_{\mb{U}},\Omega^1_{\mb{U}},\T_{\mb{U}})$ in a canonical way.
Under this identification, the isomorphism of sheaves of extended Lie 
superalgebroids induced by $\Phi$ is given by 
$\bvarphi^*:(C^\infty_{\mb{W}},\Omega^1_{\mb{W}},\T_{\mb{W}})\rightarrow
\varphi_*(C^\infty_{\mb{U}},\Omega^1_{\mb{U}},\T_{\mb{U}})$.

(b) The quotient map $\mc{V}_1|_U\rightarrow\T_{\mb{U}}$ is split as 
a morphism of sheaves of $\bb{C}$-vector spaces, and $\mc{V}|_U$ is freely 
generated by any associated sheaf of vertex superalgebroids.
Moreover, $\Phi$ is induced by an isomorphism of sheaves of vertex 
superalgebroids.
\end{lemma}

\begin{proof}
(a) At weight zero, $(\varphi,\Phi)$ defines an isomorphism of ringed spaces
$(U,C^\infty_{\mb{U}})\rightarrow(W,C^\infty_{\mb{W}})$, which is the same as
a diffeomorphism $\bvarphi:\mb{U}\rightarrow\mb{W}$.
Let $(C^\infty_{\mb{M}},\Omega,\T)$ be the presheaf of extended Lie 
superalgebroids associated to $\mc{V}$.
The following isomorphisms, induced respectively by $\Phi$ and $\bvarphi$
\beqq
\xymatrix{
  \varphi_*(C^\infty_{\mb{U}},\Omega|_U,\T|_U) &
  (C^\infty_{\mb{W}},\Omega^1_{\mb{W}},\T_{\mb{W}})
    \ar[l]_{\ph{aa}\cong}\ar[r]^\cong &
  \varphi_*(C^\infty_{\mb{U}},\Omega^1_{\mb{U}},\T_{\mb{U}})
}
\eeqq
allow us to identify $(C^\infty_{\mb{U}},\Omega|_U,\T|_U)$ with 
$(C^\infty_{\mb{U}},\Omega^1_{\mb{U}},\T_{\mb{U}})$ via identity on 
$C^\infty_{\mb{U}}$.
Since any isomorphism with 
$(C^\infty_{\mb{U}},\Omega^1_{\mb{U}},\T_{\mb{U}})$ is determined by 
its first component, the above identification is independent of the choice
of $W$ and $(\varphi,\Phi)$.

(b) The statements about $\mc{V}|_U$ are true because their analogues
for $\CDO_{p|q}|_W$ are true.
The statement about $\Phi$ is then clear.
\end{proof}

\begin{lemma} \label{globalCDO.lemma2}
Let $\mc{V}$ be a sheaf of CDOs on a smooth cs-manifold 
$\mb{M}=(M,C^\infty_{\mb{M}})$.

(a) The presheaf (in fact, sheaf) of extended Lie superalgebroids 
associated to $\mc{V}$ can be identified with 
$(C^\infty_{\mb{M}},\Omega^1_{\mb{M}},\T_{\mb{M}})$ in a canonical way.

(b) The quotient map $\mc{V}_1\rightarrow\T_{\mb{M}}$ is split as 
a morphism of sheaves of $\bb{C}$-vector spaces, and $\mc{V}$ is freely 
generated by any associated sheaf of vertex superalgebroids.
\end{lemma}

\begin{proof}
Let $\mf{U}=\{U_a\}_{a\in I}$ be an open cover of $M$ such that 
$(U_a,\mc{V}|_{U_a})$ admit isomorphisms as postulated in Definition
\ref{sheafCDO.defn}.
For $A\subset M$, let ``$A\cap\mf{U}$'' denote the open cover
$\{A\cap U_a\}_{a\in I}$ of $A$.

(a) Let $(C^\infty_{\mb{M}},\Omega,\T)$ be the \emph{presheaf} of extended
Lie superalgebroids associated to $\mc{V}$ and $U\subset M$ an arbitrary
open set.
Consider the diagram (natural in $U$)
\beqq
\xymatrix{
  &
  \Omega(U)\ar[r]^{\varepsilon\ph{aaa}}\ar@{.>}[d]_\iota &
  \check C^0(U\cap\mf{U},\Omega)\ar[r]^\delta\ar[d]^\cong &
  \check C^1(U\cap\mf{U},\Omega)\ar[d]^\cong\ar[d]^\cong \\
  0\ar[r] &
  \Omega^1_{\mb{M}}(U)\ar[r] &
  \check C^0(U\cap\mf{U},\Omega^1_{\mb{M}})\ar[r]^\delta &
  \check C^1(U\cap\mf{U},\Omega^1_{\mb{M}})
}
\eeqq
where $(\check C^*(\,\cdot\,),\delta)$ denote Cech complexes and 
the isomorphisms are given by Lemma \ref{globalCDO.lemma1}a.
By the exactness of the bottom row, the dotted arrow $\iota$ can be filled
in in a unique way.
By construction, $\iota$ is compatible with the derivations 
$C^\infty_{\mb{M}}\rightarrow\Omega$,\, 
$C^\infty_{\mb{M}}\rightarrow\Omega^1_{\mb{M}}$, and this implies $\iota$ is 
surjective.
On the other hand, since $\Omega$ is a subpresheaf of a sheaf, 
$\varepsilon$ is injective, and so is $\iota$.
Hence we have an isomorphism $\Omega\cong\Omega^1_{\mb{M}}$.
Now $\T$ must also be a sheaf.
This is a formal consequence of: 
(i) $\T(U):=\mc{V}_1(U)/\Omega(U)$ for open sets $U\subset M$,
(ii) $\mc{V}_1$ is a sheaf, and
(iii) $\Omega\cong\Omega^1_{\mb{M}}$ is a fine sheaf.
Then a diagram similar to the one above produces an isomorphism 
$\T\cong\T_{\mb{M}}$.
By construction, the isomorphisms $\Omega\cong\Omega^1_{\mb{M}}$ and
$\T\cong\T_{\mb{M}}$ respect the extended Lie superalgebroid structures.

(b) Let $\pi$ denote the quotient map $\mc{V}_1\rightarrow\T_{\mb{M}}$.
By Lemma \ref{globalCDO.lemma1}b, the restriction of $\pi$ to each $U_a$
has a splitting $s_a:\T_{\mb{M}}|_{U_a}\rightarrow\mc{V}_1|_{U_a}$.
Let $\{f_a\}_{a\in I}$ be a smooth partition of unity on $M$ subordinate 
to $\mf{U}$.
Use the operation ${}_{(-1)}:C^\infty_M\times\mc{V}_1\rightarrow\mc{V}_1$ to
define such morphisms of sheaves $(f_a)_{(-1)}s_a$ that extend from $U_a$
to $M$.
Since the said operation induces via $\pi$ the usual 
$C^\infty_{\mb{M}}$-multiplication on $\T_{\mb{M}}$, the sum
\beqq
s:=\sum_{a\in I}(f_a)_{(-1)}s_a:\T_{\mb{M}}\rightarrow\mc{V}_1
\eeqq
splits $\pi$.
Such a splitting yields a sheaf of vertex superalgebroids.

Given a sheaf of vertex superalgebroids associated to $\mc{V}$, 
its sections freely generate a \emph{presheaf} of vertex superalgebras
$\mc{V}'$.
Moreover, there is a canonical morphism of presheaves of vertex 
superalgebras $\kappa:\mc{V}'\rightarrow\mc{V}$.
Since $\kappa|_{U_a}$ are isomorphisms, so is $\kappa$ if and only if 
$\mc{V}'$ is a sheaf.
Now each weight component $\mc{V}'_k$, $k\geq 1$, admits a filtration 
whose associated graded presheaf is a sheaf 
(see \S\ref{VAoid.VA.filtration}).
It follows formally from this fact that $\mc{V}'$ is indeed a sheaf as 
desired.
\end{proof}

{\it Preparation.}
Suppose $\mb{M}$ is a smooth cs-manifold and $\nabla$ a connection on
$T\mb{M}$.
Given $X\in\T(\mb{M})$, let $\tnabla X$ denote the section of
$\t{End}\,T\mb{M}$ defined by $(\tnabla X)(Y)=\nabla_X Y-[X,Y]$ for
$Y\in\T(\mb{M})$.
Notice that if $\nabla$ is torsion-free, then $\tnabla X=\nabla X$.

\begin{theorem} \label{thm.globalCDO}
Let $\mb{M}=(M,C^\infty_{\mb{M}})$ be a smooth cs-manifold.

(a) Suppose $\nabla$ is a connection on $T\mb{M}$ with curvature operator 
$R$, and $H$ is an even $3$-form on $\mb{M}$ with $dH=\Str(R\wedge R)$.
Given such data, a sheaf of vertex superalgebroids
\beqq 
\big(C^\infty_{\mb{M}},\Omega^1_{\mb{M}},\T_{\mb{M}},
  \ast,\{\;\},\{\;\}_\Omega\big)
\eeqq
can be defined on $M$ using the following formulae
\begin{align*} 
f\ast X\;\;\; &= -(\nabla df)(X) \\
\{X,Y\}\;\; &= -\Str(\tnabla X\cdot\tnabla Y) \\
\{X,Y\}_\Omega &=
  \Str\bigg(-\nabla(\tnabla X)\cdot\tnabla Y
    +\tnabla X\cdot\iota_Y R
    -\iota_X R\cdot\tnabla Y\bigg)
  +\frac{1}{2}\iota_X\iota_Y H
\end{align*}
and it freely generats a sheaf of CDOs on $\mb{M}$, denoted by 
$\CDO_{\mb{M},\nabla,H}$.
Up to isomorphism, this construction yields all sheaves of CDOs on 
$\mb{M}$.

(b) Conformal structures on $\CDO_{\mb{M},\nabla,H}$ are in one-to-one
correspondence with even $1$-forms $\omega$ on $\mb{M}$ satisfying
$d\omega=\Str R$.
This correspondence is independent of the choice of $H$.
Given $\omega$ as described, the corresponding conformal structure, 
denoted by $\nu^\omega$, is characterized by
\beqq
L_1^\omega X
:=\nu^\omega_{(2)}X
=\Str\tnabla X-\omega(X)
\eeqq
for vector fields $X$ on $\mb{M}$.
\end{theorem}

\begin{proof}  
(a) Suppose $\mc{V}$ is a sheaf of CDOs on $\mb{M}$.
By Lemma \ref{globalCDO.lemma2}, $\mc{V}$ is freely generated by a sheaf 
of vertex superalgebroids
$\big(C^\infty_{\mb{M}},\Omega^1_{\mb{M}},\T_{\mb{M}},\ast,
\{\;\},\{\;\}_\Omega\big)$.
Let $\mf{U}=\{U_a\}_{a\in I}$ be an open cover of $M$ and 
\beqq
(\varphi_a,\Phi_a):\big(U_a,\mc{V}|_{U_a}\big)\rightarrow
  \big(W_a,\CDO_{p|q}|_{W_a}\big),\quad
W_a\subset\bb{R}^p\t{ open},\quad
a\in I
\eeqq
be isomorphisms as postulated in Definition \ref{sheafCDO.defn}.
Also let $\mb{U}_a=\mb{M}|_{U_a}$ and $\mb{W}_a=(\bb{R}^{p|q})|_{W_a}$.
By Lemma \ref{globalCDO.lemma1}, there are diffeomorphisms
$\bvarphi_a:\mb{U}_a\rightarrow\mb{W}_a$ such that $\Phi_a$ are induced by 
isomorphisms of sheaves of vertex superalgebroids of the form
\beq \label{local.VSAoid.iso}
(\bvarphi_a^*,\bvarphi_a^*\Delta_a):
\big(C^\infty_{\mb{W}_a},\Omega^1_{\mb{W}_a},\T_{\mb{W}_a},
  \ast^c,\{\;\}^c,\{\;\}^c_\Omega\big)
\rightarrow
\varphi_{a*}\big(C^\infty_{\mb{U}_a},\Omega^1_{\mb{U}_a},\T_{\mb{U}_a},
  \ast,\{\;\},\{\;\}_\Omega\big)
\eeq
where $\Delta_a:\T_{\mb{W}_a}\rightarrow\Omega^1_{\mb{W}_a}$ are some
even morphisms of sheaves on $W_a$.

Somewhat abusing notations, we will write $\varphi_a$, $\bvarphi_a$, $\Phi_a$, 
etc.~also for their restrictions to various open subsets.
For $a,a'\in I$, let $W_{aa'}=\varphi_a(U_a\cap U_{a'})$,
$\mb{W}_{aa'}=(\bb{R}^{p|q})|_{W_{aa'}}$ and
\beqq
&\varphi_{a'a}=\varphi_{a'}\circ\varphi_a^{-1}:W_{aa'}\rightarrow W_{a'a},\qquad
\bvarphi_{a'a}=\bvarphi_{a'}\circ\bvarphi_a^{-1}:
  \mb{W}_{aa'}\rightarrow\mb{W}_{a'a} & \\
&(\varphi_{a'a},\Phi_{a'a})
=(\varphi_{a'},\Phi_{a'})\circ(\varphi_a,\Phi_a)^{-1}:
\big(W_{aa'},\CDO_{p|q}|_{W_{aa'}}\big)\rightarrow
\big(W_{a'a},\CDO_{p|q}|_{W_{a'a}}\big) &
\eeqq
Recall the notations in \S\ref{MC.WZ} and write $g_{\bvarphi_{a'a}}$,
$\theta_{\bvarphi_{a'a}}$, $WZ_{\bvarphi_{a'a}}$ more simply as
$g_{a'a}$, $\theta_{a'a}$, $WZ_{a'a}$. \vss
According to \S\ref{sec.CDO.coorchange}, 
$\Phi_{a'a}=(\bvarphi_{a'a})^*_{\xi_{a'a}}$ for some unique even $2$-forms
$\xi_{a'a}$ on $\mb{W}_{aa'}$ with $d\xi_{a'a}=WZ_{a'a}$, and it is 
induced by an isomorphism of sheaves of vertex superalgebroids
$(\bvarphi_{a'a}^*,\Delta_{a'a})$, where
\beqq
\Delta_{a'a}=\Delta_{\bvarphi_{a'a},\xi_{a'a}}:
\T_{\mb{W}_{a'a}}\rightarrow(\varphi_{a'a})_*\Omega^1_{\mb{W}_{aa'}}
\eeqq
is defined as in Theorem \ref{CDO.iso}.
The definition of $\Phi_{a'a}$ given above is equivalent to
\begin{align}
(\bvarphi_{a'a}^*,\Delta_{a'a})
&=(\varphi_{a'a})_*(\bvarphi_a^*,\bvarphi_a^*\Delta_a)^{-1}\circ
  (\bvarphi_{a'}^*,\bvarphi_{a'}^*\Delta_{a'}) \nonumber \\
\Leftrightarrow\qquad\qquad
\Delta_{a'a}&=\bvarphi_{a'a}^*\circ\Delta_{a'}
  -(\varphi_{a'a})_*\Delta_a\circ\bvarphi_{a'a}^* \label{resolve.Deltaxi}
\end{align}
For $a,a',a''\in I$, let $W_{aa'a''}=\varphi_a(U_a\cap U_{a'}\cap U_{a''})$
and $\mb{W}_{aa'a''}=(\bb{R}^{p|q})|_{W_{aa'a''}}$.
In $W_{aa'a''}$ we have
\beqq
\bvarphi_{a''a}=\bvarphi_{a''a'}\circ\bvarphi_{a'a},\qquad
(\bvarphi_{a''a})^*_{\xi_{a''a}}
  =(\varphi_{a''a'})_*(\bvarphi_{a'a})^*_{\xi_{a'a}}
  \circ(\bvarphi_{a''a'})^*_{\xi_{a''a'}}
\eeqq
According to \S\ref{sec.CDO.coorchange}, the latter is equivalent to
\beq \label{xi.cocycle}
\xi_{a''a'}=\xi_{a'a}+\bvarphi_{a'a}^*\xi_{a''a'}+\sigma_{a''a'a}
\eeq
where $\sigma_{a''a'a}=\sigma_{\bvarphi_{a''a'},\,\bvarphi_{a'a}}\in
\Omega^2(\mb{W}_{aa'a''})$ is defined as in Theorem \ref{CDO.iso}.

\begin{lemma} \label{resolve.Deltaxi.soln}
Given $\bvarphi_{a'a}$ and $\xi_{a'a}$ for $a,a'\in I$ as above (which
determine $\Delta_{a'a}$), a collection of even morphisms of sheaves 
$\Delta_a:\T_{\mb{W}_a}\rightarrow\Omega^1_{\mb{W}_a}$ satisfy 
(\ref{resolve.Deltaxi}) if and only if they are of the form
\beqq
\Delta_a(X)
=\e_i\e_{ij}\e_j^{1+|X|}(\d_j X^i)(\Gamma_a)^j_{\ph{i}i}
+\frac{1}{2}\iota_X\Str(\Gamma_a\otimes\Gamma_a)
+\frac{1}{2}\iota_X B_a
+O_a(X)
\eeqq
for homogeneous $X$, where: \vss \\
\indent\;$\cdot$\;
$\Gamma_a\in\Omega^1(\mb{W}_a)\otimes\mf{gl}(p|q)$ are even, 
i.e.~$|(\Gamma_a)^i_{\ph{i}j}|=|b^i|+|b^j|$, and 
\beq \label{Gamma.soln}
g_{a'a}^{-1}\cdot\bvarphi_{a'a}^*\Gamma_{a'}\cdot g_{a'a}-\Gamma_a
=-\theta_{a'a}
\eeq
\indent\;$\cdot$\;
$B_a\in\Omega^2(\mb{W}_a)$ are even and
\beq \label{B.soln}
\bvarphi_{a'a}^*B_{a'}-B_a
=-\xi_{a'a}
-\Str(\theta_{a'a}\wedge\Gamma_a)
\eeq
and $O_a:\T_{\mb{W}_a}\rightarrow\Omega^1_{\mb{W}_a}$ are even and
$\bvarphi_{a'a}^*\circ O_{a'}=(\varphi_{a'a})_*O_a\circ\bvarphi_{a'a}^*$.
\end{lemma}

\begin{proof}
If we assume $\Delta_a$ are first-order differential operators, we may 
write
\beqq
\Delta_a(X)
=\e_i\e_{ij}\e_j^{1+|X|}(\d_j X^i)(\Gamma_a)^j_{\ph{i}i}
+\frac{1}{2}\iota_X(S_a+B_a)
\eeqq
for some $\mf{gl}(p|q)$-valued $1$-forms $\Gamma_a$,
symmetric $(0,2)$-tensors $S_a$ and $2$-forms $B_a$ on $\mb{W}_a$;
their parities are dictated by that of $\Delta_a$.
Plugging this into (\ref{resolve.Deltaxi}), namely
\beqq
\bvarphi_{a'a}^*\Delta_{a'}(X)-\Delta_a(\bvarphi_{a'a}^*X)=\Delta_{a'a}(X)
\eeqq
results in three sets of equations: 
(\ref{Gamma.soln}), (\ref{B.soln}) and
\beqq
\bvarphi_{a'a}^*S_{a'}-S_a=
-\Str(\Gamma_a\otimes\theta_{a'a})
-\Str(\theta_{a'a}\otimes\Gamma_a)
+\Str(\theta_{a'a}\otimes\theta_{a'a}).
\eeqq
By (\ref{Gamma.soln}), the last set of equations is satisfied by 
$S_a=\Str(\Gamma_a\otimes\Gamma_a)$.
Observe that once we have a solution to (\ref{resolve.Deltaxi}), any other 
solution differs precisely by a term $O_a$ with the properties stated in
the lemma.
\end{proof}

\noindent
{\it Proof of Theorem \ref{thm.globalCDO} continued.}
Consider the formula of $\Delta_a$ obtained in Lemma 
\ref{resolve.Deltaxi.soln}.
The condition on the term $O_a$ lets us define a map
$O:\T_{\mb{M}}\rightarrow\Omega^1_{\mb{M}}$ such that
$O(\bvarphi_a^*X)=\bvarphi_a^*O_a(X)$ for $a\in I$.
By Lemma \ref{lemma.newVAoid}, $O$ determines an isomorphism of sheaves
of vertex superalgebroids 
\beqq
(\t{id},-O):
\big(C^\infty_{\mb{M}},\Omega^1_{\mb{M}},\T_{\mb{M}},
  \ast,\{\;\},\{\;\}_\Omega\big)
\rightarrow
\big(C^\infty_{\mb{M}},\Omega^1_{\mb{M}},\T_{\mb{M}},
  \ast',\{\;\}',\{\;\}'_\Omega\big)
\eeqq
whose composition with (\ref{local.VSAoid.iso}) equals
\beqq
\varphi_{a*}(\t{id},-O)|_{U_a}\circ(\bvarphi_a^*,\bvarphi_a^*\Delta_a)
=(\bvarphi_a^*,\bvarphi_a^*\Delta_a-\bvarphi_a^*O_a).
\eeqq
Therefore up to isomorphism of sheaves of CDOs, we may assume $O_a=0$.
The following lemma concerns the other ingredients in the formula of
$\Delta_a$.

\begin{lemma} \label{xi.B.H.equiv}
Assume that $U_a$, $a\in I$, are contractible.
Given $\bvarphi_{a'a}$ for $a,a'\in I$ as above, the existence of 
the following are equivalent: \vspace{-0.06in}
\begin{itemize}
\item[(i)] 
$\xi_{a'a}\in\Omega^2(\mb{W}_{a'a})$ that are even, and satisfy 
$d\xi_{a'a}=WZ_{a'a}$ and (\ref{xi.cocycle}) \vspace{-0.06in}
\item[(ii)]
$\Gamma_a\in\Omega^1(\mb{W}_a)\otimes\mf{gl}(p|q)$ and 
$B_a\in\Omega^2(\mb{W}_a)$ that are even, and satisfy (\ref{Gamma.soln}) 
and
\beq \label{dB.Cech}
\bvarphi_{a'a}^*dB_{a'}-dB_a=-\bvarphi_{a'a}^*CS(\Gamma_{a'})+CS(\Gamma_a)
\eeq
where $CS(\Gamma_a)\in\Omega^3(\mb{W}_a)$ is defined below \vspace{-0.06in}
\item[(iii)]
a connection $\nabla$ on $T\mb{M}$ and $H\in\Omega^3(\mb{M})$ that is even
and satisfies
\beq \label{dH}
dH=\Str(R\wedge R)
\eeq
where $R$ is the curvature operator of $\nabla$
\end{itemize}
\end{lemma}

\begin{proof}
First, a collection of $\Gamma_a\in\Omega^1(\mb{W}_a)\otimes\mf{gl}(p|q)$ 
that are even and satisfy (\ref{Gamma.soln}) is equivalent to a connection
$\nabla$ on $T\mb{M}$. 
Indeed, the two are related via
\beq \label{conn.Gamma}
\nabla(\bvarphi_a^*\d_i)
=\e_i\e_{ij}\,\bvarphi_a^*\big((\Gamma_a)^j_{\ph{i}i}\otimes\d_j\big)
\eeq
for $i=1,\ldots,p+q$.
The curvature operator $R\in\Omega^2(\mb{M},\t{End}\,T\mb{M})$ of $\nabla$
is locally given by
\beqq
R(\bvarphi_a^*\d_i)
=\e_i\e_{ij}\,\bvarphi_a^*\big((R_a)^j_{\ph{i}i}\otimes\d_j\big),\qquad
R_a=d\Gamma_a+\Gamma_a\wedge\Gamma_a
\eeqq
whose tensoriality means
$g_{a'a}^{-1}\cdot\bvarphi_{a'a}^*R_{a'}\cdot g_{a'a}=R_a$.
Define the following even $3$-forms on $\mb{W}_a$
\beqq
CS(\Gamma_a):=\Str(\Gamma_a\wedge R_a)
  -\frac{1}{3}\Str(\Gamma_a\wedge\Gamma_a\wedge\Gamma_a).
\eeqq
Notice that $d\,CS(\Gamma_a)=\Str(R_a\wedge R_a)$ and (\ref{Gamma.soln}) 
implies
\beq \label{CS.Cech}
\bvarphi_{a'a}^*CS(\Gamma_{a'})-CS(\Gamma_a)
  =WZ_{a'a}+d\,\Str(\theta_{a'a}\wedge\Gamma_a).
\eeq
Now we prove the equivalences.

(i) $\Rightarrow$ (ii):
Choose a connection $\nabla$ on $T\mb{M}$ and define $\Gamma_a$ as in 
(\ref{conn.Gamma}).
Then $\Gamma_a$ satisfy (\ref{Gamma.soln}).
The right hand side of (\ref{B.soln}), after being pulled back by 
$\bvarphi_a^*$, defines a $1$-cochain in the \v Cech complex 
$\check C^*(\mf{U},\Omega^2_{\mb{M}})$;
it is a cocycle by (\ref{xi.cocycle}) and (\ref{Gamma.soln}).
Since $\check C^*(\mf{U},\Omega^2_{\mb{M}})$ is acyclic, we may choose 
such even $2$-forms $B_a$ on $\mb{W}_a$ that satisfy (\ref{B.soln}).
Then (\ref{dB.Cech}) follows from $d\xi_{a'a}=WZ_{a'a}$ and 
(\ref{CS.Cech}).

(ii) $\Rightarrow$ (i):
Define $\xi_{a'a}$ using (\ref{B.soln}).
Then (\ref{Gamma.soln}) implies (\ref{xi.cocycle}).
On the other hand, (\ref{dB.Cech}) and (\ref{CS.Cech}) together imply
$d\xi_{a'a}=WZ_{a'a}$.

(ii) $\Rightarrow$ (iii):
Define $\nabla$ as in (\ref{conn.Gamma}).
By (\ref{dB.Cech}), there is a global even $3$-form $H$ on $\mb{M}$ with 
\beq \label{H.dB.CS}
H|_{U_a}=\bvarphi_a^*\big(dB_a+CS(\Gamma_a)\big).
\eeq
Then (\ref{dH}) follows from $d\,CS(\Gamma_a)=\Str(R_a\wedge R_a)$.

(iii) $\Rightarrow$ (ii):
Define $\Gamma_a$ as in (\ref{conn.Gamma}).
Then $\Gamma_a$ satisfy (\ref{Gamma.soln}).
The $3$-forms $H|_{U_a}-\bvarphi_a^*CS(\Gamma_a)$ are closed because of 
(\ref{dH}) and the fact that $d\,CS(\Gamma_a)=\Str(R_a\wedge R_a)$.
Since $U_a$ are contractible, we may choose such even $2$-forms $B_a$ on 
$\mb{W}_a$ that satisfy (\ref{H.dB.CS}), which implies (\ref{dB.Cech}).
\end{proof}

\noindent
{\it Proof of Theorem \ref{thm.globalCDO} continued.}
Now compute the maps $\ast^\dag$, $\{\;\}_0^\dag$, $\{\;\}_1^\dag$.
In view of (\ref{local.VSAoid.iso}), the restrictions of the three maps
to $U_a$ are given by
\begin{align}
f\ast X\;\;\; &= \bvarphi_a^*\Big(
  f_a\ast^c X_a+\Delta_a(f_a X_a)-f_a\Delta_a(X_a)\Big) \nonumber \\
\{X,Y\}\;\; &= \bvarphi_a^*\Big(
  \{X_a,Y_a\}^c-\Delta_a(X_a)(Y_a)-(-1)^{|X||Y|}\Delta_a(Y_a)(X_a)\Big)
  \label{globalCDO.VAoid456.Delta} \\
\{X,Y\}_\Omega &= \bvarphi_a^*\Big(
  \{X_a,Y_a\}^c_\Omega
  -L_{X_a}\Delta_a(Y_a)
  +(-1)^{|X||Y|}L_{Y_a}\Delta_a(X_a)
  -d\,\Delta_a(X_a)(Y_a)+\Delta_a([X_a,Y_a])\Big) \nonumber
\end{align}
where $f_a,X_a,Y_a$ are such that $f=\bvarphi_a^*f_a$, $X=\bvarphi_a^*X_a$, 
$Y=\bvarphi_a^*Y_a$.
To evaluate (\ref{globalCDO.VAoid456.Delta}), apply the formulae of 
$\ast^c$, $\{\;\}^c$, $\{\;\}^c_\Omega$ in (\ref{CDO.VAoid456}), and
that of $\Delta_a$ in Lemma \ref{resolve.Deltaxi.soln} (with $O_a=0$).
Also use the data $\Gamma_a$, $B_a$ in the formula of $\Delta_a$ to
define a connection $\nabla$ on $T\mb{M}$ and an even $3$-form $H$ on
$\mb{M}$ as in (\ref{conn.Gamma}) and (\ref{H.dB.CS});
by the proof of Lemma \ref{xi.B.H.equiv} they satisfy (\ref{dH}).
A lengthy but straightforward computation then yields the formulae of
$\ast$, $\{\;\}$, $\{\;\}_\Omega$ stated in the theorem.
This proves the last statement of part (a).

It remains to argue that the construction described in the theorem always
produces a sheaf of CDOs on $\mb{M}$.
Notations in this paragraph will have the same meaning as above.
Choose a covering $\mf{U}=\{U_a\}_{a\in I}$ of $M$ by contractible open
sets, and diffeomorphisms $\bvarphi_a:\mb{U}_a\rightarrow\mb{W}_a$;
let $\bvarphi_{a'a}=\bvarphi_{a'}\circ\bvarphi_a^{-1}$.
Starting with the given data $\nabla$, $H$, define $\Gamma_a$, $B_a$ as in
the proof of Lemma \ref{xi.B.H.equiv}, and then $\Delta_a$ as in Lemma 
\ref{resolve.Deltaxi.soln} (with $O_a=0$).
By the same computation mentioned before, $\Delta_a$ and the given formulae
of $\ast$, $\{\;\}$, $\{\;\}_\Omega$ satisfy 
(\ref{globalCDO.VAoid456.Delta}).
Then by Lemma \ref{lemma.newVAoid}, $\ast$, $\{\;\}$, $\{\;\}_\Omega$
define a sheaf of vertex superalgebroids equipped with the isomorphisms
(\ref{local.VSAoid.iso}).
Its freely generated sheaf of vertex superalgebras is therefore a sheaf 
of CDOs.

(b) Use the notations in (a).
Suppose $\nu$ is a conformal structure on $\mc{V}$.
For $a\in I$
\beqq
\nu|_{U_a}=\Phi_a\big(\nu^{\omega_a}\big)
\eeqq
for some even closed $1$-forms $\omega_a$ on $\mb{W}_a$.
For $a,a'\in I$, the isomorphism $\Phi_{a'a}=(\bvarphi_{a'a})^*_{\xi_{a'a}}$
sends $\nu^{\omega_{a'}}$ to $\nu^{\omega_a}$.
By Proposition \ref{conformal.coorchange}, this is equivalent to 
the relation
\beqq
\bvarphi_{a'a}^*\omega_{a'}-\omega_a=\Str\theta_{a'a}
=-\bvarphi_{a'a}^*\Str\Gamma_{a'}+\Str\Gamma_a
\eeqq
where the second equality is given by (\ref{Gamma.soln}).
Hence there is an even $1$-form $\omega$ on $\mb{M}$ with 
\beqq
\omega|_{U_a}=\bvarphi_a^*(\omega_a+\Str\Gamma_a).
\eeqq
Since $d\omega_a=0$ and $d\,\Str\Gamma_a=\Str R_a$, we have 
$d\omega=\Str R$.
Observe that the construction of $\omega$ from $\nu$ is reversible.
To relate $\nu$ and $\omega$ more explicitly, we compute
$\Phi_a(\nu^{\omega_a})$ as follows:
\begin{align}
\nu|_{U_a}
& =\e_i\big(\bvarphi_a^*\d_i+\bvarphi_a^*\Delta_a(\d_i)\big)_{-1}
    (\bvarphi_a^*db^i)
  +\frac{1}{2}(\bvarphi_a^*\omega_a)_{-2}\bf{1} 
  \nonumber \\
& =\e_i\left(\frac{\d}{\d\bvarphi_a^i}\right)_{-1}d\bvarphi_a^i
  +\frac{1}{2}\Str\bigg((\bvarphi_a^*\Gamma_a)_{-1}(\bvarphi_a^*\Gamma_a)
    -(\bvarphi_a^*\Gamma_a)_{-2}\mb{1}\bigg)
  +\frac{1}{2}\omega_{-2}\mb{1}
  \label{globalconformal.coor}
\end{align}
where we first recall that $\Phi_a$ is induced by (\ref{local.VSAoid.iso})
and then use Lemma \ref{resolve.Deltaxi.soln}.
The computation does not depend on $B_a$, hence not on $H$.
Let $L_n^\omega=\nu_{(n+1)}$.
Using (\ref{globalconformal.coor}) we have
\begin{align*}
L_1^\omega X\big|_{U_a}
& =(d\bvarphi_a^i)_1\left(\frac{\d}{\d\bvarphi_a^i}\right)_0 X
  +\Str(\bvarphi_a^*\Gamma_a)_1 X
  -\omega_1 X \\
& =d\bvarphi_a^i\left(\left[\frac{\d}{\d\bvarphi_a^i},X\right]\right)
  +\Str(\bvarphi_a^*\Gamma_a)(X)
  -\omega(X)
\end{align*}
for vector fields $X$.
The sum of the first two terms is a local expression for $\Str\tnabla X$.
\end{proof}

{\it Remarks.}
(i) Given an open set $U\subset M$, let $\mb{U}=\mb{M}|_U$.
The vertex superalgebra $\CDO_{\mb{M},\nabla,H}(U)$ will also be written as
$\CDO_{\nabla,H}(\mb{U})$.
A conformal structure $\nu$ on $\CDO_{\mb{M},\nabla,H}$ restricts to
a conformal structure $\nu|_U$ on $\CDO_{\nabla,H}(\mb{U})$ of
central charge $2(p-q)$.
(ii) By definition, there are canonical identifications
\beqq
(\CDO_{\mb{M},\nabla,H})_0=C^\infty_{\mb{M}},\qquad
(\CDO_{\mb{M},\nabla,H})_1=\Omega^1_{\mb{M}}\oplus\T_{\mb{M}}.
\eeqq
Consider the following $\bb{C}$-bilinear operation for each $k\geq 0$
\beqq
C^\infty_{\mb{M}}\times(\CDO_{\mb{M},\nabla,H})_k\rightarrow
  (\CDO_{\mb{M},\nabla,H})_k,\qquad
(f,v)\mapsto f_0 v=f_{(-1)}v.
\eeqq
For $k>0$, this operation does not make $(\CDO_{\mb{M},\nabla,H})_k$ 
a $C^\infty_{\mb{M}}$-module,
\footnote{
For example, we have $f_0 X=fX+f\ast X=fX-(\nabla df)(X)$ for vector
fields $X$.
}
but it induces a $C^\infty_{\mb{M}}$-module structure on an associated
graded sheaf $gr(\CDO_{\mb{M},\nabla,H})_k$.
Given a $C^\infty_{\mb{M}}$-module $\mc{E}$, we use the notation
$\widehat{\t{Sym}}_t\mc{E}$ for the formal sum
$\sum_{n=0}^\infty t^n\widehat{\t{Sym}}{}^n\mc{E}$,
where $t$ is  a formal variable and $\widehat{\t{Sym}}{}^n\mc{E}$
is the $n$-fold graded symmetric tensor power of
$\mc{E}$ over $C^\infty_{\mb{M}}$.
There is an isomorphism of $C^\infty_{\mb{M}}$-modules
\beqq
C^\infty_{\mb{M}}\oplus
\bigoplus_{k=1}^\infty q^k gr(\CDO_{\mb{M},\nabla,H})_k
\cong\bigotimes_{\ell=1}^\infty
  \widehat{\t{Sym}}_{q^\ell}(\Omega^1_{\mb{M}}\oplus\T_{\mb{M}}).
\eeqq
For more details of the vertex superalgebra structure of 
$\CDO_{\mb{M},\nabla,H}$, consult \S\ref{VAoid.VA} and 
\S\ref{VAoid.VA.filtration}.

\begin{theorem} \label{thm.globalCDO.iso}
Let $\mb{M}$ be a smooth cs-manifold, $\nabla$ a connection on $T\mb{M}$ 
with curvature operator $R$, and $H,H'$ even $3$-forms on $\mb{M}$ with
$dH=dH'=\Str(R\wedge R)$.
Define $\ast$, $\{\;\}$, $\{\;\}_\Omega$ (resp.~$\{\;\}'_\Omega$)
using $\nabla$ and $H$ (resp.~$H'$) as in Theorem \ref{thm.globalCDO}a.

(a) There is a one-to-one correspondence:
\beqq
\left\{\hspace{-0.05in}
\begin{array}{c}
  \t{isomorphisms of sheaves of CDOs }
  \CDO_{\mb{M},\nabla,H}\rightarrow\CDO_{\mb{M},\nabla,H'} \vss \\ 
  \t{whose weight-zero components are identity on }C^\infty_{\mb{M}}
\end{array}\hspace{-0.05in}
\right\}
\stackrel{\sim}{\longleftrightarrow}
\left\{\hspace{-0.05in}
\begin{array}{c}
  B\in\Omega^2(\mb{M}),\t{even},\\
  dB=H'-H
\end{array}\hspace{-0.05in}
\right\}
\eeqq
Given $B$ as above, the corresponding isomorphism, denoted by 
$\mathrm{id}_B$, is induced by an isomorphism between the associated 
sheaves of vertex superalgebroids
\beqq
(\mathrm{id},\Delta_B):
\big(C^\infty_{\mb{M}},\Omega^1_{\mb{M}},\T_{\mb{M}},
  \ast,\{\;\},\{\;\}_\Omega\big)
\rightarrow
\big(C^\infty_{\mb{M}},\Omega^1_{\mb{M}},\T_{\mb{M}},
  \ast,\{\;\},\{\;\}'_\Omega\big)
\eeqq
where the map $\Delta_B:\T_{\mb{M}}\rightarrow\Omega^1_{\mb{M}}$ is 
given by $\Delta_B(X)=\frac{1}{2}\iota_X B$.

(b) The isomorphism $\mathrm{id}_B$ preserves the correspondence in
Theorem \ref{thm.globalCDO}b, i.e.~$\mathrm{id}_B(\nu^\omega)=\nu^\omega$.
\end{theorem}

\begin{proof}
(a) If an isomorphism between the two sheaves of CDOs equals the identity 
on $C^\infty_{\mb{M}}$, then it induces the identity on the sheaf of extended
Lie superalgebroids $(C^\infty_{\mb{M}},\Omega^1_{\mb{M}},\T_{\mb{M}})$, and
is therefore determined by an isomorphism of sheaves of vertex 
superalgebroids of the form
\beqq
(\mathrm{id},\Delta):
\big(C^\infty_{\mb{M}},\Omega^1_{\mb{M}},\T_{\mb{M}},
  \ast,\{\;\},\{\;\}_\Omega\big)
\rightarrow
\big(C^\infty_{\mb{M}},\Omega^1_{\mb{M}},\T_{\mb{M}},
  \ast,\{\;\},\{\;\}'_\Omega\big).
\eeqq
By definition, the even map
$\Delta:\T_{\mb{M}}\rightarrow\Omega^1_{\mb{M}}$
has to satisfy precisely the following equations:
\beqq
&\Delta(fX)=f\Delta(X),\qquad
\Delta(Y)(X)=-(-1)^{|X||Y|}\Delta(X)(Y)& \\
&\ds L_X\Delta(Y)-(-1)^{|X||Y|}L_Y\Delta(X)+d\,\Delta(X)(Y)-\Delta([X,Y])
=\{X,Y\}_\Omega-\{X,Y\}'_\Omega&
\eeqq
According to the first two equations, $B(X,Y):=2\Delta(X)(Y)$ defines 
an even $2$-form $B$ on $\mb{M}$.
Then the last equation can be rewritten as
\beqq
\iota_X\iota_Y dB=\iota_X\iota_Y(H'-H).
\eeqq

(b) Since the said correspondence is independent of $H$, this is clear.
This also follows from the local expression (\ref{globalconformal.coor}) of
$\nu^\omega$ and the graded symmetry of $(d\bvarphi_a^i)_{-1}d\bvarphi_a^j$.
\end{proof}

\begin{example} \label{vbcs.CDO}
{\bf Sheaves of CDOs on $\Pi E$.}
Let $M$ be a smooth manifold, 
$E\rightarrow M$ a smooth $\bb{C}$-vector bundle and
$\mb{M}=\Pi E$ as a smooth cs-manifold.
The canonical pullback embeds $\Omega^*(M)$ into
$\Omega^*(\mb{M})$ quasi-isomorphically.~\cite{QFS.susy}
Choose connections $\nabla^M$ on $TM$ and $\nabla^E$ on $E$,
which determine a connection $\nabla$ on $T\mb{M}$
in the sense of \S\ref{sec.cs.conn};
denote by $R^M$, $R^E$ and $R$ the corresponding curvature
tensors.
As stated in Lemma \ref{lemma.str}, we have
\beqq
\Str(R\wedge R)=\Tr(R^M\wedge R^M)-\Tr(R^E\wedge R^E),\qquad
\Str R=\Tr R^M-\Tr R^E.
\eeqq
By Theorems \ref{thm.globalCDO}a and \ref{thm.globalCDO.iso}a,
$\mb{M}$ admits sheaves of CDOs if and only if
$p_1(TM)-ch_2(E)$ vanishes in de Rham cohomology,
and their isomorphism classes form
an $H^3(M;\bb{C})$-torsor.
By Theorem \ref{thm.globalCDO}b, the sheaves of CDOs possess
conformal structures if and only if $c_1(E)$ vanishes
in de Rham cohomology as well.
\end{example}

\begin{example} \label{cdR}
{\bf The smooth chiral de Rham complex.}
Consider the case $E=TM\otimes\bb{C}$ in 
the previous example.
Both obstructions are now trivial, so that $\mb{M}$ always
admits sheaves of CDOs equipped with conformal structures.
In particular, we may define a sheaf of CDOs
$\CDO_{\mb{M},\nabla}=\CDO_{\mb{M},\nabla,0}$ using
the trivial $3$-form and a conformal structure $\nu=\nu^0$
using the trivial $1$-form.

Let $J$ and $Q$ be the vector fields on $\mb{M}$
defined in \S\ref{sec.cs.vf} and Example \ref{dRcs}.
Regarded as elements of $\CDO_\nabla(\mb{M})$ of
weight $1$, they satisfy
\beqq
2Q_0^2
=[Q_0,Q_0]
=[Q,Q]_0+(\{Q,Q\}_\Omega)_0
=0,\qquad
[J_0,Q_0]
=[J,Q]_0+(\{J,Q\}_\Omega)_0
=Q_0.
\eeqq
In view of the formulae in Theorem \ref{thm.globalCDO}a and
Lemma \ref{lemma.1form.0mode}, the two equations follow from
(\ref{Q.bracket}) and Lemma \ref{lemma.str.Q}, with
the second also requiring Lemmas \ref{lemma.ops.J}b
and \ref{lemma.str}c.
\footnote{
In fact, we are assuming that $\nabla^M$ is Levi-Civita.
The torsion-free condition is used to obtain various
formulae in Example \ref{dRcs} and subsequently
Lemma \ref{lemma.str.Q}, while orthogonality ensures that
the right hand side of Lemma \ref{lemma.str}c vanishes.
}
Moreover, we have
\beqq
Q_0\nu=-\frac{1}{2}T^2(L_1 Q)=0,\qquad
J_0\nu=-\frac{1}{2}T^2(L_1 J)=0.
\eeqq
In view of Theorem \ref{thm.globalCDO}b, the two
equations follow from Lemmas \ref{lemma.str.Q}
and \ref{lemma.ops.J}a respectively.
Therefore, with $J_0$ as the grading operator and $Q_0$
as the differential, $\CDO_\nabla(\mb{M})$ becomes
a differential graded conformal vertex superalgebra.
\footnote{
For a description of a richer structure on
$\CDO_\nabla(\mb{M})$, see \cite{BHS}.
}
\end{example}

\begin{lemma} \label{lemma.1form.0mode}
Consider a sheaf of CDOs $\CDO_{\mb{M},\nabla,H}$ on a smooth cs-manifold
$\mb{M}$ constructed as in Theorem \ref{thm.globalCDO}.
Given $\alpha\in\Omega^1(\mb{M})$, we have $\alpha_0=0$ on
$\CDO_{\nabla,H}(\mb{M})$ if and only if $d\alpha=0$.
\end{lemma}

\begin{proof}
Since $\alpha_0$ is a derivation, it acts trivially on the entire vertex
superalgebra $\CDO_{\nabla,H}(\mb{M})$ if and only if it acts trivially on
functions and vector fields.
For $f\in C^\infty(\mb{M})$, we always have $\alpha_0 f=0$.
For $X\in\T(\mb{M})$, we compute
\beqq
\alpha_0 X=\pm[X_{-1},\alpha_0]\mb{1}=\pm(L_X\alpha-d\iota_X\alpha)
  =\pm\iota_X d\alpha 
\eeqq
which proves the assertion.
\end{proof}

\newpage

\setcounter{equation}{0}
\section{Chiral Dolbeault Algebras}

Applying the description of CDOs obtained in
Theorem \ref{thm.globalCDO}, we study a vertex algebraic
analogue of the Dolbeault complex of a complex manifold.
This provides a new point of view on the relation between
CDOs and elliptic genera.

\begin{subsec} \label{sec.Dolbcs}
{\bf Dolbeault cs-manifolds.}
Let $M$ be a complex manifold, $TM$ its holomorphic tangent
bundle, $E$ a holomorphic vector bundle over $M$, and
$\mb{M}=\Pi(\overline{TM}\oplus E)$ as a smooth cs-manifold.
Let $d=\dim_{\bb{C}}M$ and $r=\t{rank}\,E$.
Under the identification
\beq \label{Dolb.fcn}
C^\infty(\mb{M})\cong\Omega^{0,*}(M;\wedge^*E^\vee)
\eeq
vector fields on $\mb{M}$ correspond to derivations of
the $(0,*)$-forms on $M$ valued in $\wedge^*E^\vee$.
In particular, let
\beq \label{Dolb.JJQ}
\left\{\begin{array}{l}
  J^r \\
  J^\ell \\
  Q 
\end{array}\right\}
=
\begin{array}{c}
 \t{the vector field on }\mb{M} \\
 \t{corresponding to the }
\end{array}
\left\{\begin{array}{l}
  \t{Dolbeault degree} \\
  \t{exterior degree in }\wedge^*E^\vee \\
  \t{Dolbeault operator }\dbar\otimes 1
\end{array}\right\}
\eeq
For more discussion of $\T(\mb{M})$, see Example \ref{Dolbcs}.
\end{subsec}

\begin{subsec} \label{sec.DolbCDO}
{\bf Sheaves of CDOs on $\mb{M}$.}
Choose connections $\nabla^M$ on $TM$ and $\nabla^E$ on $E$
such that both are of type $(1,0)$ and $\nabla^M$ is
torsion-free (see footnote \ref{torsionfree}).
Let $\nabla$ be the induced connection on $T\mb{M}$
defined as in Example \ref{Dolbcs}.
Denote by $R^M$, $R^E$ and $R$ the respective curvature
tensors.
Notice that the canonical pullback embeds $\Omega^*(M)$ into
$\Omega^*(\mb{M})$ quasi-isomorphically \cite{QFS.susy} and
recall Lemma \ref{lemma.str.Dolb}.

Assume that $ch_2(TM)-ch_2(E)=0$ in de Rham cohomology and
choose $H\in\Omega^3(M)$ such that
\beq \label{DolbCDO.H}
dH=\Str(R\wedge R)
  =\Tr(R^M\wedge R^M)-\Tr(R^E\wedge R^E).
\eeq
By Theorems \ref{thm.globalCDO}a and
\ref{thm.globalCDO.iso}a, this determines a sheaf of CDOs
$\CDO_{\mb{M},\nabla,H}$ and every sheaf of CDOs on $\mb{M}$
is up to isomorphism of this form.
Assume also that $c_1(TM)-c_1(E)=0$ in de Rham cohomology
and choose $\omega\in\Omega^1(M)$ such that
\beq \label{DolbCDO.omega}
d\omega=\Str R=\Tr R^M-\Tr R^E.
\eeq
By Theorem \ref{thm.globalCDO}b, this determines a conformal
structure $\nu^\omega$ on $\CDO_{\mb{M},\nabla,H}$ of
central charge $2(d-r)$.
\end{subsec}

\begin{theorem} \label{thm.DolbQ}
Regard $Q$ as an element of $\CDO_{\nabla,H}(\mb{M})$ of
weight $1$.
The odd derivation $Q_0$:

(a) is a differential if and only if $H$ has no $(1,2)$-
or $(0,3)$-part, and

(b) respects the conformal structure $\nu^\omega$ if and
only if $\omega$ has no $(0,1)$-part.
\end{theorem}

\begin{proof}
(a) The supercommutator of $Q_0$ with itself is given by
\beqq
2Q_0^2
=[Q_0,Q_0]
=[Q,Q]_0+(\{Q,Q\}_\Omega)_0
=\frac{1}{2}(\iota_Q\iota_Q H)_0
\eeqq
where the last step follows from (\ref{DolbQ.bracket}),
Theorem \ref{thm.globalCDO}a, Lemma \ref{lemma.str.DolbQ}
and Lemma \ref{lemma.1form.0mode}.
By Lemma \ref{lemma.1form.0mode} again, $Q_0^2$ vanishes if and
only if $\iota_Q\iota_Q H$ is closed.
In view of the identity
\beqq
2\iota_Q\iota_Q H
=L_{J^r}\iota_Q\iota_Q H
=\iota_{J^r}d\iota_Q\iota_Q H
\eeqq
$\iota_Q\iota_Q H$ can only be closed when it is in fact
trivial.
When applied to a differential form on $M$, $\iota_Q\iota_Q$
picks out those components of type $(i,j)$ with $j\geq 2$.

(b) Applying $Q_0$ to $\nu^\omega$ yields
\beqq
Q_0(\nu^\omega)
=-[\nu^\omega_{(-1)},Q_0]\mb{1}
=-\frac{1}{2}T^2 L^\omega_1Q
=\frac{1}{2}T^2\omega(Q)
\eeqq
where the last step follows from Theorem \ref{thm.globalCDO}b
and Lemma \ref{lemma.str.DolbQ}.
Hence $Q_0$ annihilates $\nu^\omega$ if and only if
$T\omega(Q)=d\omega(Q)=0$.
In view of the identity
\beqq
\omega(Q)=J^r\omega(Q)=\iota_{J^r}d\omega(Q)
\eeqq
$\omega(Q)$ can only be constant when it is in fact trivial.
When applied to a differential form on $M$, $\iota_Q$ picks
out those components of type $(i,j)$ with $j\geq 1$.
\end{proof}

\begin{defn} \label{hol.Chern}
For each $n\geq 0$, let $\Omega^{n,\t{cl}}_{M,\t{hol}}$
(resp.~$\Omega^{n,\t{cl}}_M$) denote the sheaf of holomorphic
(resp.~smooth) closed $n$-forms on $M$ and define an element
\beqq
c^{\t{hol}}_n(E)\in H^n(\Omega^{n,\t{cl}}_{M,\t{hol}})
\eeqq
as follows.
Since $\nabla^E$ is of type $(1,0)$, its curvature $R^E$ has
only $(2,0)$- and $(1,1)$-parts.
Thus the $n$-th Chern form $c_n(\nabla^E)$ lives in
$\Omega^{n+\ast,\ast}(M)$.
Consider the diagram of fine resolutions of sheaves:
\beqq
\xymatrix{
  0\ar[r] &
  \Omega^{n,\t{cl}}_{M,\t{hol}}\ar[r]\ar[d]^\cap &
  (\Omega^{n+*,*}_M,d)\ar[d]^\cap \\
  0\ar[r] &
  \Omega^{n,\t{cl}}_M\ar[r] &
  (\Omega^{n+*}_M,d)
}
\eeqq
In light of this diagram, $c_n(\nabla^E)$ represents an element 
``$c^{\t{hol}}_n(E)$'' of $H^n(\Omega^{n,\t{cl}}_{M,\t{hol}})$,
whose image under
\beq \label{hol.sm.cl.chlgy}
H^n(\Omega^{n,\t{cl}}_{M,\t{hol}})
\rightarrow H^n(\Omega^{n,\t{cl}}_M)
\cong H^{2n}(M;\bb{C})
\eeq
is the $n$-th Chern class $c_n(E)$.
More generally, if $C(E)$ is a polynomial in the Chern classes
$c_n(E)$, denote by $C^{\t{hol}}(E)$ the corresponding polynomial
in $c_n^{\t{hol}}(E)$.
The following result relates some of these cohomology classes to
the conditions encountered in Theorem \ref{thm.DolbQ}.
\end{defn}

\begin{prop} \label{prop.obstructions}
There exists:

(a) $H\in\Omega^3(M)$ satisfying (\ref{DolbCDO.H}) and
$H^{1,2}=H^{0,3}=0$ if and only if
$ch^{\mathrm{hol}}_2(TM)-ch^{\mathrm{hol}}_2(E)=0$;

(b) $\omega\in\Omega^1(M)$ satisfying (\ref{DolbCDO.omega}) and
$\omega^{0,1}=0$ if and only if
$c^{\mathrm{hol}}_1(TM)-c^{\mathrm{hol}}_1(E)=0$.
\end{prop}

\begin{proof}
Recall Definition \ref{hol.Chern}.
Statement (a) holds because the said element of 
$H^2(\Omega^{2,\t{cl}}_{M,\t{hol}})$ is represented, via the fine
resolution
\beqq
\xymatrix{
  0\ar[r] &
  \Omega^{2,\t{cl}}_{M,\t{hol}}\ar[r] &
  \Omega^{2,0}_M\ar[r]^{d\ph{aaaa}} &
  \Omega^{3,0}_M\oplus\Omega^{2,1}_M\ar[r]^{d\ph{aaaa}} &
  \Omega^{4,0}_M\oplus\Omega^{3,1}_M\oplus\Omega^{2,2}_M
    \ar[r]^{\ph{aaaaaaa}d} &
  \cdots 
}
\eeqq
by the right hand side of (\ref{DolbCDO.H}) up to a constant
factor.
Statement (b) is similar.
\end{proof}

{\it Remark.}
In the case $M$ is K\"ahler, (\ref{hol.sm.cl.chlgy}) is injective,
as it can be identified with the inclusion
\beqq
\bigoplus_{p\geq 0, p+q=n}\mc{H}^{n+p,q}\hookrightarrow\mc{H}^{2n}
\eeqq
where $\mc{H}^{n+p,q}$, $\mc{H}^{2n}$ are the spaces of harmonic
$(n+p,q)$- and $2n$-forms respectively.
Thus the conditions in Proposition \ref{prop.obstructions} become
equivalent to $ch_2(TM)-ch_2(E)=0$ and $c_1(TM)-c_1(E)=0$.

\begin{subsec} \label{sec.fermion}
{\bf Fermion numbers.}
The eigenvalues of $J^r_0$ and $J^\ell_0$ on
$\CDO_{\mb{M},\nabla,H}$ will be referred to respectively
as right (i.e.~antiholomorphic) and left (i.e.~holomorphic)
fermion numbers.
Recall from (\ref{Dolb.JJQ}) that in weight $0$,
these numbers correspond to the exterior degrees in
$\wedge^*\overline{TM}{}^\vee$ and $\wedge^*E^\vee$.
\end{subsec}

\begin{prop} \label{prop.Q.Jr}
The operator $Q_0$ always increases right fermion numbers
by $1$ if and only if the line bundle $\det TM$ is flat.
\end{prop}

\begin{proof}
Denote by $\bar\nabla^M$ the connection on $\overline{TM}$
induced by $\nabla^M$ and $\bar R^M$ its curvature tensor.
The commutator between $J^r_0$ and $Q_0$ is given by
\beqq
[J^r_0,Q_0]
=[J^r,Q]_0+(\{J^r,Q\}_\Omega)_0
=Q_0-(\iota_Q\Tr\bar R^M)_0
\eeqq
which follows from (\ref{DolbQ.bracket}), Theorem
\ref{thm.globalCDO}a, Lemmas \ref{lemma.Jr}b,
\ref{lemma.str.Dolb}c and \ref{lemma.str.DolbQ}.
Hence by Lemma \ref{lemma.1form.0mode}, $Q_0$ is compatible
with right fermion numbers if and only if
$\iota_Q\Tr\bar R^M$ is closed.
By the same argument used in the proof of Theorem
\ref{thm.DolbQ}, $\iota_Q\Tr\bar R^M$ can only be closed
when it is in fact trivial.
Since $\nabla^M$ is of type $(1,0)$, $\bar R^M$ has only
$(1,1)$- and $(0,2)$-parts, so that $\iota_Q\Tr\bar R^M=0$
if and only if $\Tr\bar R^M=0$.
\end{proof}

\begin{prop} \label{prop.Q.Jl}
The operator $Q_0$ respects left fermion numbers
if and only if $\Tr R^E$ has no $(1,1)$-part.
\end{prop}

\begin{proof}
The commutator between $J^\ell_0$ and $Q_0$ is given by
\beqq
[J^\ell_0,Q_0]
=[J^\ell,Q]_0+(\{J^\ell,Q\}_\Omega)_0
=-(\iota_Q\Tr R^E)_0
\eeqq
which follows from (\ref{DolbQ.bracket}), Theorem
\ref{thm.globalCDO}a, Lemmas \ref{lemma.Jl}b,
\ref{lemma.str.Dolb}c and \ref{lemma.str.DolbQ}.
Hence by Lemma \ref{lemma.1form.0mode}, $Q_0$ commutes
with $J^\ell_0$ if and only if $\iota_Q\Tr R^E$ is closed.
By the same argument used in the proof of Theorem
\ref{thm.DolbQ}, $\iota_Q\Tr R^E$ can only be closed
when it is in fact trivial.
Since $\nabla^E$ is of type $(1,0)$, $R^E$ has only
$(2,0)$- and $(1,1)$-parts, so that $\iota_Q$ picks out
the $(1,1)$-part.
\end{proof}

{\it Remark.}
Given a hermitian metric on $E$, there exists a unique
unitary connection $\nabla^E$ of type $(1,0)$, and its
curvature $R^E$ is of pure type $(1,1)$.~\cite{Wells}
If $\CDO_{\nabla,H}(\mb{M})$ has been defined using this
$\nabla^E$, then $Q_0$ respects left fermion numbers if and
only if $\Tr R^E=0$, i.e.~the line bundle $\det E$ is flat.

\begin{corollary}
Suppose $Q_0^2=0$ holds, so that $(\CDO_{\nabla,H}(\mb{M}),Q_0)$
is a differential vertex superalgebra.

(a) If $\det TM\cong\det E$ as holomorphic line bundles,
$(\CDO_{\nabla,H}(\mb{M}),Q_0)$ is a differential conformal
vertex superalgebra.

(b) If $\det TM$ is flat, the grading by right fermion numbers
makes $(\CDO_{\nabla,H}(\mb{M}),Q_0)$ a differential
graded vertex superalgebra.

(c) If $\det E$ is flat, left fermion numbers are well-defined
on the cohomology of $(\CDO_{\nabla,H}(\mb{M}),Q_0)$.
\end{corollary}

\begin{proof}
(a) Under the assumption we may compare the induced connections
$\det\nabla^M$ and $\det\nabla^E$ via the isomorphism, and they
differ by a $(1,0)$-form $\omega$.
This implies (\ref{DolbCDO.omega}) and, by
Theorem \ref{thm.DolbQ}b, $Q_0(\nu^\omega)=0$.
(b)-(c) These are simply restatements of Propositions
\ref{prop.Q.Jr} and \ref{prop.Q.Jl}.
\end{proof}

For the rest of this section, $M$ is always compact.

\begin{theorem} \label{thm.cDolb.chlgy}
Suppose $Q_0^2=0$ holds and consider the vertex superalgebra
\beqq 
V=H\big(\CDO_{\nabla,H}(\mb{M}),Q_0\big).
\eeqq
Let $q$ be a formal variable.
There is an identity of formal power series
\beq \label{part.func.L0}
\mathrm{Str}_V(q^{L_0})
=\int_M\mathrm{Td}(TM)\,ch\left(
  \bigotimes_{n=1}^\infty\mathrm{Sym}_{q^n}(TM\oplus TM^\vee)\otimes
  \bigotimes_{n=1}^\infty\wedge_{-q^n}E\otimes
  \bigotimes_{n=0}^\infty\wedge_{-q^n}E^\vee
\right).
\eeq
Let $y$ be another formal variable.
If $\det E$ is flat, there is a more refined identity
\beqq
\mathrm{Str}_V(y^{J^\ell_0}q^{L_0})
=\int_M\mathrm{Td}(TM)\,ch\left(
  \bigotimes_{n=1}^\infty\mathrm{Sym}_{q^n}(TM\oplus TM^\vee)\otimes
  \bigotimes_{n=1}^\infty\wedge_{-y^{-1}q^n}E\otimes
  \bigotimes_{n=0}^\infty\wedge_{-yq^n}E^\vee
\right).
\eeqq
\end{theorem}

\begin{proof}
By Proposition \ref{prop.Q.Jl}, if $\det R^E$ is flat, $J^\ell_0$
is well-defined on $V$.
Otherwise, set $y=1$ whenever it appears in the proof below.

Observe that $Q_0$ respects the filtration on
$\CDO_{\nabla,H}(\mb{M})$ described in
\S\ref{VAoid.VA.filtration} and induces the operator $L_Q$ on
the associated graded space $gr(\CDO_{\nabla,H}(\mb{M}))$.
Let
\beqq
V'=H\big(gr(\CDO_{\nabla,H}(\mb{M})),L_Q\big).
\eeqq
The quantity we want to compute can be rephrased as follows:
\begin{align}
&\t{ supertrace of }y^{J^l_0}q^{L_0}\t{ on }V \nonumber \\
=&\t{ supertrace of }y^{J^l_0}q^{L_0}\t{ on }V' \nonumber \\
=&\t{ supertrace of }y^{J^l_0}\t{ on }
  H\left(\bigotimes_{n=1}^\infty
  \widehat{\t{Sym}}_{q^n}\big(\Omega^1(\mb{M})\oplus\T(\mb{M})\big),
  L_Q\right) \label{cDolb.chlgy.char}
\end{align}
where the graded symmetric tensor products are taken over
$C^\infty(\mb{M})$.
Recall the local coordinates
\beqq
(\t{Re}z^1,\t{Im}z^1,\cdots,\t{Re}z^d,\t{Im}z^d,
  \bar\zeta^1,\cdots,\bar\zeta^d,
  \ve^1,\cdots,\ve^r)
\eeqq
defined in Example \ref{Dolbcs}.
To compute (\ref{cDolb.chlgy.char}), consider the following
subspaces of $\Omega^1(\mb{M})$ and $\T(\mb{M})$:
\begin{align*}
\Omega^{1,0}(\mb{M})\ph{'}&=\big\{
  \alpha\in\Omega^1(\mb{M})\t{ s.t. }
  \alpha\t{ is locally a }C^\infty\t{-linear combination of }
  dz^i, d\ve^k\big\} \\
\Omega^{1,0}(\mb{M})'&=\big\{
  \alpha\in\Omega^1(\mb{M})\t{ s.t. }
  \alpha\t{ is locally a }C^\infty\t{-linear combination of }
  dz^i\big\} \\
\T^{1,0}(\mb{M})\ph{'}&=\big\{
  X\in\T(\mb{M})\t{ s.t. }
  X\t{ is locally a }C^\infty\t{-linear combination of }
  \d/\d z^i, \d/\d\ve^k\big\} \\
\T^{1,0}(\mb{M})'&=\big\{
  X\in\T(\mb{M})\t{ s.t. }
  X\t{ is locally a }C^\infty\t{-linear combination of }
  \d/\d\ve^k\big\}
\end{align*}

\begin{lemma} \label{1forms.10forms.qiso}
The following inclusions
\beqq
\big(\Omega^{1,0}(\mb{M}),L_Q\big)\hookrightarrow
  \big(\Omega^1(\mb{M}),L_Q\big),
\qquad
\big(\T^{1,0}(\mb{M}),L_Q\big)\hookrightarrow
  \big(\T(\mb{M}),L_Q\big)
\eeqq
are quasi-isomorphisms.
\end{lemma}

\begin{proof}
Denote both of the projections
$\Omega^1(\mb{M})\rightarrow\Omega^{1,0}(\mb{M})$ and
$\T(\mb{M})\rightarrow\T^{1,0}(\mb{M})$ by $\pi^{1,0}$.
It suffices to show that $\t{id}-\pi^{1,0}$ are null
homotopic.
Define $G:\Omega^1(\mb{M})\rightarrow\Omega^1(\mb{M})$
and $G:\T(\mb{M})\rightarrow\T(\mb{M})$ locally by
\beqq
G\alpha
=(-1)^{|\alpha|}
  \alpha\Big(\frac{\d}{\d\bar\zeta^i}\Big)d\bar z^i,\qquad
GX
=(-1)^{|X|}
  d\bar z^i(X)\frac{\d}{\d\bar\zeta^i}
\eeqq
and notice that the expressions are independent of
local coordinates.
By a calculation we have
\beqq
L_Q G+GL_Q=\t{id}-\pi^{1,0}
\eeqq
on both $\Omega^1(\mb{M})$ and $\T(\mb{M})$, as desired.
\end{proof}

\begin{lemma} \label{10forms.filtration}
There are natural filtrations on $(\Omega^{1,0}(\mb{M}),L_Q)$
and $(\T^{1,0}(\mb{M}),L_Q)$ whose associated graded complexes
are isomorphic respectively to
\beqq
\big(\Omega^{0,*}(M;E'),\dbar\big)
\quad\t{and}\quad
\big(\Omega^{0,*}(M;E''),\dbar\big)
\eeqq
where $E'=\wedge^*E\otimes(TM^\vee\oplus E^\vee)$ and
$E''=\wedge^*E\otimes(TM\oplus E)$.
\end{lemma}

\begin{proof}
There are identifications defined by the following
local expressions
\beqq
\begin{array}{ll}
\Omega^{1,0}(\mb{M})'\cong
  C^\infty(\mb{M})\otimes_{C^\infty(M)}\Omega^{1,0}(M),&
dz^i\mapsto 1\otimes dz^i \vss \\
\Omega^{1,0}(\mb{M})/\Omega^{1,0}(\mb{M})'
  \cong C^\infty(\mb{M})\otimes_{C^\infty(M)}\Gamma(E^\vee),\quad& 
d\ve^k\t{ mod }\Omega^{1,0}(\mb{M})'\mapsto 1\otimes\ve^k \vss \\
\T^{1,0}(\mb{M})'\cong
  C^\infty(\mb{M})\otimes_{C^\infty(M)}\Gamma(E),&
\d/\d\ve^k\mapsto 1\otimes\ve_k \vss \\
\T^{1,0}(\mb{M})/\T^{1,0}(\mb{M})'\cong
  C^\infty(\mb{M})\otimes_{C^\infty(M)}\T^{1,0}(M),&
\d/\d z^i\t{ mod }\T^{1,0}(\mb{M})'\mapsto 1\otimes\d/\d z^i
\end{array}
\eeqq
and $C^\infty(\mb{M})$-linearity.
Then it remains to recall (\ref{Dolb.fcn}).
\end{proof}

\noindent
{\it Proof of Theorem \ref{thm.cDolb.chlgy} continued.}
Now we apply Lemmas \ref{1forms.10forms.qiso} and
\ref{10forms.filtration} to compute (\ref{cDolb.chlgy.char}).
Since (graded symmetric) tensor products of quasi-isomorphic
complexes are quasi-isomorphic and filtrations on complexes
induce filtrations on their (graded symmetric) tensor products,
we have
\begin{align*}
(\ref{cDolb.chlgy.char})=&
\t{ supertrace of }y^{J^l_0}\t{ on }
  H\left(\bigotimes_{n=1}^\infty
  \widehat{\t{Sym}}_{q^n}\big(\Omega^{1,0}(\mb{M})\oplus
  \T^{1,0}(\mb{M})\big),L_Q\right) \\
=&\t{ supertrace of }y^{J^l_0}\t{ on }
  H\left(\bigotimes_{n=1}^\infty
  \widehat{\t{Sym}}_{q^n}\Omega^{0,*}(M;E'\oplus E''),
  \dbar\right) \\
=&\t{ sdim}\;H\left(\Omega^{0,*}
  \bigg(M;\wedge_{-y}E^\vee\otimes
  \bigotimes_{n=1}^\infty\widehat{\t{Sym}}_{q^n}
  \big(TM\oplus TM^\vee\oplus(-y^{-1})E\oplus(-y)E^\vee\big)\bigg),
  \dbar\right) \\
=&\t{ sdim}\;H\left(\Omega^{0,*}\bigg(M;
  \bigotimes_{n=1}^\infty\t{Sym}_{q^n}(TM\oplus TM^\vee)\otimes
  \bigotimes_{n=1}^\infty\wedge_{-y^{-1}q^n}E\otimes
  \bigotimes_{n=0}^\infty\wedge_{-yq^n}E^\vee\bigg),\dbar\right)
\end{align*}
Notice that the graded symmetric tensor products in the second expression
are taken over $\Omega^{0,*}(M;\wedge^*E^\vee)$.
To finish the proof, apply the Hirzebruch-Riemann-Roch Theorem.
\end{proof}

{\it Remark.}
In terms of the Chern roots $x_1,\ldots,x_d$ of $TM$ and
$x^E_1,\ldots,x^E_r$ of $E$, we may write the integrand in
(\ref{part.func.L0}) as follows
{\footnotesize
\begin{align*}
&\prod_{i=1}^d\left(
  \frac{x_i}{1-e^{-x_i}}
  \prod_{n=1}^\infty\frac{1}{(1-q^n e^{x_i})(1-q^n e^{-x_i})}
\right)\cdot
\prod_{j=1}^r\left(
  (1-e^{-x^E_j})
  \prod_{n=1}^\infty(1-q^n e^{x^E_j})(1-q^n e^{-x^E_j})
\right) \\
=\;&
\prod_{i=1}^d\left(
  \frac{x_i/2}{\sinh(x_i/2)}
  \prod_{n=1}^\infty\frac{1}{(1-q^n e^{x_i})(1-q^n e^{-x_i})}
\right)\cdot
\prod_{j=1}^r\left(
  2\sinh\frac{x^E_j}{2}
  \prod_{n=1}^\infty(1-q^n e^{x^E_j})(1-q^n e^{-x^E_j})
\right)\cdot
\frac{e^{\frac{1}{2}c_1(TM)}}{e^{\frac{1}{2}c_1(E)}}
\end{align*}}

\noindent
If $c_1(TM)=c_1(E)$, this expression lives in $H^{4*}(M;\bb{C})$
if $r$ is even, or in $H^{4*+2}(M;\bb{C})$ if $r$ is odd, so that
$\Str_V(q^{L_0})=0$ whenever $d+r$ is odd.

\begin{example} \label{E=0}
{\bf The case $E=0$.}
By Theorem \ref{thm.DolbQ}a and Proposition
\ref{prop.obstructions}a, there exists a differential
vertex superalgebra
\beqq
\big(\CDO_{\nabla,H}(\mb{M}),Q_0\big),\qquad
\mb{M}=\Pi\,\overline{TM}
\eeqq
if and only if $ch_2^{\t{hol}}(TM)=0$;
denote its cohomology by $V$.
By Theorem \ref{thm.cDolb.chlgy}
\beqq
\t{Str}_V(q^{L_0})
=\int_M e^{\frac{1}{2}c_1(TM)}\cdot W(TM_{\bb{R}})\cdot 
  \bigg(\prod_{n=1}^\infty\frac{1}{1-q^n}\bigg)^{2d}
\eeqq
where $W(TM_{\bb{R}})$ is the Witten class of the real
tangent bundle of $M$.
By Theorem \ref{thm.DolbQ}b and Proposition
\ref{prop.obstructions}b, if $c_1^{\t{hol}}(TM)=0$ as well,
$V$ is conformal with central charge $2d$.
Then, writing $q=e^{2\pi i\tau}$, we have
\beq \label{cDolb.Witten}
\t{char}\,V
=q^{-d/12}\,\t{Str}_V(q^{L_0})
=\frac{W(M)}{\Delta(\tau)^{d/12}}
\eeq
where $W(M)$ is the Witten genus of $M$ and
\beqq
\Delta(\tau)=q\prod_{n=1}^\infty(1-q^n)^{24}.
\eeqq
The condition $c_1(TM)=c_2(TM)=0$ guarantees that $W(M)$
is a modular form of weight $d$, while $\Delta(\tau)$ is
a modular form of weight $12$, both over $SL(2,\bb{Z})$.
The expression in (\ref{cDolb.Witten}) is the conjectured
$S^1$-equivariant index of the Dirac operator on the free loop space
$LM$.~\cite{Witten.index}
\end{example}

\begin{example} \label{E=TM}
{\bf The case $E=TM$.}
By Theorem \ref{thm.DolbQ} and Proposition
\ref{prop.obstructions}, there always exists a differential
conformal vertex superalgebra
\beqq
\big(\CDO_{\nabla,H}(\mb{M}),Q_0\big),\qquad
\mb{M}=\Pi(\overline{TM}\oplus TM)
\eeqq
with no central charge;
denote its cohomology by $V$.
By Theorem \ref{thm.cDolb.chlgy}, we have
\beqq
\t{Str}_V(q^{L_0})=\chi(M)
\eeqq
and, if $\det TM$ is flat, also have
\beqq
y^{-d}\,\t{Str}_V(y^{J^l_0}q^{L_0})=Ell_{y,q}(M)
\eeqq
namely the two-variable elliptic genus of $M$.~\cite{BL}
In particular, writing $q=e^{2\pi i\tau}$, we have
the special value
\beq \label{cDolb.elliptic}
\t{Str}_V((-1)^{J^l_0}q^{L_0})
  =\frac{\t{Och}(M)}{\epsilon(\tau)^{d/4}}
\eeq
where $\t{Och}(M)$ is the Ochanine elliptic genus of $M$ and 
\beqq
\epsilon(\tau)
=\frac{1}{16}\prod_{n=1}^\infty\bigg(\frac{1-q^n}{1+q^n}\bigg)^8
\eeqq
respectively a modular form of weight $d$ and a modular form of
weight $4$ over $\Gamma_0(2)\subset SL(2,\bb{Z})$.
The expression in (\ref{cDolb.elliptic}) is the $S^1$-equivariant
signature of $LM$.~\cite{HBJ}
\end{example}

\begin{example} \label{E=l}
{\bf The case $E=\det TM$.}
Let $c=c_1(TM)$ and $c^{\t{hol}}=c_1^{\t{hol}}(TM)$.
By Theorem \ref{thm.DolbQ}a and Proposition
\ref{prop.obstructions}a, there exists a differential
vertex superalgebra
\beqq
\big(\CDO_{\nabla,H}(\mb{M}),Q_0\big),\qquad
\mb{M}=\Pi(\overline{TM}\oplus\det TM)
\eeqq
if and only if
\beq \label{hol.stringc.odd}
ch_2^{\t{hol}}(TM)-\frac{1}{2}(c^{\t{hol}})^2=0;
\eeq
denote its cohomology by $V$.
By Theorem \ref{thm.DolbQ}b and Proposition
\ref{prop.obstructions}b, $V$ is always conformal with
central charge $2(d-1)$.
By Theorem \ref{thm.cDolb.chlgy} and the remark below
its proof, we have
\beq \label{Wittenc.odd}
\t{Str}_V(q^{L_0})
=2\int_M
  W(TM_{\bb{R}})\,\sinh\frac{c}{2}\,
  \prod_{n=1}^\infty\frac{(1-q^n e^c)(1-q^n e^{-c})}{(1-q^n)^2}
  \cdot\bigg(\prod_{n=1}^\infty\frac{1}{1-q^n}\bigg)^{2(d-1)}
\eeq
which always vanishes if $d$ is even.
Now assume $d$ is odd.
This case provides a geometric interpretation of
the notions introduced in \cite{CHZ} for certain spin$^c$
manifolds of (real) dimension $2\t{ mod }4$.
Firstly, condition (\ref{hol.stringc.odd}) implies that
$M$ is \emph{rationally string$^c$} in the sense of
{\it loc.~cit.}, namely
\beq \label{stringc.odd}
ch_2(TM)-\frac{1}{2}c^2=0\quad\t{in }H^4(M;\bb{C}).
\eeq
In the case $M$ is K\"ahler, (\ref{hol.stringc.odd}) and
(\ref{stringc.odd}) are equivalent, as remarked after
the proof of Proposition \ref{prop.obstructions}.
Secondly, writing $q=e^{2\pi i\tau}$, we have
\beqq
\t{char}\,V
=q^{-(d-1)/12}\t{Str}_V(q^{L_0})
=\frac{2W_c(M)}{\Delta(\tau)^{(d-1)/12}}
\eeqq
where $W_c(M)$ is the \emph{generalized Witten genus} of $M$
defined in {\it loc.~cit.}
\footnote{
To recover the expression for $W_c(M)$ in \cite{CHZ}, notice
that they write $q=e^{\pi i\tau}$, and the factor $\sinh(c/2)$ in
(\ref{Wittenc.odd}) may be replaced by $e^{c/2}$ since $d=\dim_{\bb{C}}M$
is odd.
}
The string$^c$ condition (\ref{stringc.odd}) guarantees that
$W_c(M)$ is a modular form of weight $d-1$ over $SL(2,\bb{Z})$.
\end{example}

\begin{example} \label{E=l2-l}
{\bf The case $E=(\det TM)^{\otimes 2}-\det TM$.}
\footnote{
The results obtained above may be \emph{formally} applied
to a virtual holomorphic vector bundle $E=E_1-E_2$.
This amounts to using
``$C^\infty(\mb{M})$''$:=\Gamma(\wedge^*E_1^\vee
\otimes\t{Sym}^*E_2^\vee)$
in the construction of ``$\CDO_{\nabla,H}(\mb{M})$.''
}
Let $c=c_1(TM)$ and $c^{\t{hol}}=c_1^{\t{hol}}(TM)$.
By Theorem \ref{thm.DolbQ}a and Proposition
\ref{prop.obstructions}a, there exists a differential
vertex superalgebra
\beqq
\big(\CDO_{\nabla,H}(\mb{M}),Q_0\big),\qquad
\mb{M}=``\Pi\big(\overline{TM}\oplus
  ((\det TM)^{\otimes 2}-\det TM)\big)\t{''}
\eeqq
if and only if
\beq \label{hol.stringc.even}
ch_2^{\t{hol}}(TM)-\frac{3}{2}(c^{\t{hol}})^2=0;
\eeq
denote its cohomology by $V$.
By Theorem \ref{thm.DolbQ}b and Proposition
\ref{prop.obstructions}b, $V$ is always conformal
with central charge $2d$.
By Theorem \ref{thm.cDolb.chlgy} and the remark below
its proof, we have
\begin{align}
\t{Str}_V(q^{L_0})
&=2\int_M
  W(TM_{\bb{R}})\,\cosh\frac{c}{2}\,
  \prod_{n=1}^\infty
    \frac{(1-q^n e^{2c})(1-q^n e^{-2c})}{(1-q^n e^c)(1-q^n e^{-c})}
  \cdot
  \bigg(\prod_{n=1}^\infty\frac{1}{1-q^n}\bigg)^{2d} \nonumber \\
&=2\int_M
  W(TM_{\bb{R}})\,\cosh\frac{c}{2}\,
  \prod_{n=1}^\infty\bigg[
    (1-q^{n-\frac{1}{2}}e^c)(1+q^{n-\frac{1}{2}}e^c)(1+q^n e^c)
  \nonumber \\
&\qquad\qquad\qquad\cdot
    (1-q^{n-\frac{1}{2}}e^{-c})(1+q^{n-\frac{1}{2}}e^{-c})(1+q^n e^{-c})
  \bigg]\cdot
  \bigg(\prod_{n=1}^\infty\frac{1}{1-q^n}\bigg)^{2d}
\label{Wittenc.even}
\end{align}
which always vanishes if $d$ is odd.
Now assume $d$ is even.
This case provides a geometric interpretation of
the notions introduced in \cite{CHZ} for certain spin$^c$
manifolds of (real) dimension divisible by $4$.
Firstly, condition (\ref{hol.stringc.even}) implies that
$M$ is \emph{rationally string$^c$} in the sense of
{\it loc.~cit.}, namely
\beq \label{stringc.even}
ch_2(TM)-\frac{3}{2}c^2=0\quad\t{in }H^4(M;\bb{C}).
\eeq
In the case $M$ is K\"ahler, (\ref{hol.stringc.even}) and
(\ref{stringc.even}) are equivalent, as remarked after
the proof of Proposition \ref{prop.obstructions}.
Secondly, writing $q=e^{2\pi i\tau}$, we have
\beqq
\t{char}\,V
=q^{-d/12}\t{Str}_V(q^{L_0})
=\frac{2W_c(M)}{\Delta(\tau)^{d/12}}
\eeqq
where $W_c(M)$ is the \emph{generalized Witten genus} of $M$
defined in {\it loc.~cit.}
\footnote{
To recover the expression for $W_c(M)$ in \cite{CHZ}, notice that
they write $q=e^{\pi i\tau}$, and the factor $\cosh(c/2)$ in
(\ref{Wittenc.even}) may be replaced by $e^{c/2}$ since
$d=\dim_{\bb{C}}M$ is even.
}
The string$^c$ condition (\ref{stringc.even}) guarantees that
$W_c(M)$ is a modular form of weight $d$ over $SL(2,\bb{Z})$.
\end{example}

\newpage

\titleformat{\section}[block]
  {\sc\large\filcenter}
  {Appendix \S\thesection.}{.5em}{}

\appendix

\setcounter{equation}{0}
\section{Vertex Algebroids} \label{app.VAoid}

The notion of a vertex algebroid, introduced in \cite{GMS2},
captures the part of structure of a vertex algebra involving
only the two lowest weights.
In this appendix, we review the category of vertex algebroids, 
the forgetful functor from vertex algebras to vertex algebroids,
and its adjoint functor.
Some examples are given, including the construction of local
smooth CDOs.

\begin{defn}
An \emph{extended Lie algebroid} $(A,\Omega,\T)$ consists of \\
\indent $\cdot$\;
a commutative, associative $\bb{C}$-algebra with unit $(A,\mathbf 1)$ \\
\indent $\cdot$\;
two $A$-modules $\Omega$ and $\T$ \\
\indent $\cdot$\;
an $A$-derivation $d:A\rightarrow\Omega$ whose image generates $\Omega$ 
as an $A$-module \\
\indent $\cdot$\;
a Lie bracket $[\;]$ on $\T$ \\
\indent $\cdot$\;
an $A$-linear homomorphism of Lie algebras $\T\rightarrow\t{End}\,A$,
denoted $X\mapsto X$ \\
\indent $\cdot$\;
a $\bb{C}$-linear homomorphism of Lie algebras 
$\T\rightarrow\t{End}\,\Omega$, denoted $X\mapsto L_X$ \\
\indent $\cdot$\;
an $A$-bilinear pairing $\Omega\times\T\rightarrow A$, denoted 
$(\alpha,X)\mapsto\alpha(X)$ \\
Furthermore, we require that \\
\indent $\cdot$\;
the $\T$-actions on $A$ and $\Omega$ commute with $d$ \\
\indent $\cdot$\;
the $\T$-actions on $A,\Omega$ and $\T$ (via $[\;]$) satisfy the Leibniz 
rule w.r.t.~$A$-multiplication \\
\indent $\cdot$\;
$df(X)=Xf$ for $f\in A$, $X\in\T$
\end{defn}

\begin{defn}
A \emph{morphism of extended Lie algebroids} 
$\varphi:(A,\Omega,\T)\rightarrow(A',\Omega',\T')$ is a map of triples 
that respects the extended Lie algebroid structures.
Composition of morphisms is defined in the obvious way.
\end{defn}

\begin{defn}
A \emph{vertex algebroid} $(A,\Omega,\T,\ast,\{\;\},\{\;\}_\Omega)$ consists of
an extended Lie algebroid $(A,\Omega,\T)$ and three $\bb{C}$-bilinear maps
\beqq
\ast:A\times\T\rightarrow\Omega,\qquad
\{\;\}:\T\times\T\rightarrow A,\qquad
\{\;\}_\Omega:\T\times\T\rightarrow\Omega
\eeqq
that satisfy the following identities \vs \\
\indent $\cdot$\; 
$\{X,Y\}=\{Y,X\}$ \vss \\
\indent $\cdot$\; 
$d\{X,Y\}=\{X,Y\}_\Omega+\{Y,X\}_\Omega$ \vss \\
\indent $\cdot$\; 
$(fg)\ast X-f\ast(gX)-f(g\ast X)=-(Xf)dg-(Xg)df$ \vss \\
\indent $\cdot$\; 
$\{X,fY\}-f\{X,Y\}=-(f\ast Y)(X)-YXf$ \vss \\
\indent $\cdot$\;
$\{X,fY\}_\Omega-f\{X,Y\}_\Omega=-L_X(f\ast Y)+(Xf)\ast Y+f\ast[X,Y]$ \vss \\
\indent $\cdot$\;
$X\{Y,Z\}-\{[X,Y],Z\}-\{Y,[X,Z]\}=\{X,Y\}_\Omega(Z)+\{X,Z\}_\Omega(Y)$ \vss \\
\indent $\cdot$\;
$L_X\{Y,Z\}_\Omega-L_Y\{X,Z\}_\Omega+L_Z\{X,Y\}_\Omega
+\{X,[Y,Z]\}_\Omega-\{Y,[X,Z]\}_\Omega-\{[X,Y],Z\}_\Omega$ \\
\indent\qquad $=d\big(\{X,Y\}_\Omega(Z)\big)$ \vss \\
for $f,g\in A$ and $X,Y,Z\in\T$.
\end{defn}

{\it Remark.}
This definition is slightly different from but equivalent to the original
one in \cite{GMS2}.
What we denote by $\ast,\{\;\},\{\;\}_\Omega$ equal respectively 
$-\gamma,\langle\;\rangle,-c+\frac{1}{2}d\circ\langle\;\rangle$ in 
their notations.

\begin{defn} \label{VAoid.mor}
A \emph{morphism of vertex algebroids} 
\beqq
(\varphi,\Delta):(A,\Omega,\T,\ast,\{\;\},\{\;\}_\Omega)\rightarrow
(A',\Omega',\T',\ast',\{\;\}',\{\;\}'_\Omega)
\eeqq
consists of a morphism of extended Lie algebroids 
$\varphi:(A,\Omega,\T)\rightarrow(A',\Omega',\T')$ and a $\bb{C}$-linear map 
$\Delta:\T\rightarrow\Omega'$ such that \vs \\
\indent $\cdot$\;
$\varphi f\ast'\varphi X-\varphi(f\ast X)=
\Delta(fX)-(\varphi f)\Delta(X)$ \vss \\
\indent $\cdot$\;
$\{\varphi X,\varphi Y\}'-\varphi\{X,Y\}=
-\Delta(X)(\varphi Y)-\Delta(Y)(\varphi X)$ \vss \\
\indent $\cdot$\;
$\{\varphi X,\varphi Y\}'_\Omega-\varphi\{X,Y\}_\Omega=
-L_{\varphi X}\Delta(Y)+L_{\varphi Y}\Delta(X)
-d\big(\Delta(X)(\varphi Y)\big)+\Delta([X,Y])$ \vss \\
for $f\in A$ and $X,Y\in\T$.
Composition of morphisms is given by
\beqq
(\varphi',\Delta')\circ(\varphi,\Delta)=
(\varphi'\varphi,\,\varphi'\Delta+\Delta'\varphi|_\T).
\eeqq
\end{defn}

\begin{subsec}
{\bf The vertex algebroid associated to a vertex algebra 
(and a ``splitting'').}
Given a vertex algebra $(V,\mathbf 1,T,Y)$, let
\beqq
A:=V_0,\qquad
\Omega:=A_{(-1)}(TA),\qquad
\T:=V_1/\Omega.
\eeqq
Choose a splitting $s:\T\rightarrow V_1$ of the quotient map to obtain
an identification of vector spaces
\beq \label{V1.iden}
\Omega\oplus\T\cong V_1,\qquad (\alpha,X)\mapsto\alpha+s(X).
\eeq
The vertex algebra structure on $V$ involving only the two lowest weights 
consists of an element $\mathbf 1\in V_0$, 
a linear map $T:V_0\rightarrow V_1$, and eight bilinear maps
\beqq
{}_{(i+j-k-1)}:V_i\times V_j\rightarrow V_k,\quad i,j,k=0,1
\eeqq
satisfying a set of (Borcherds) identities.
These data, when rephrased in terms of the identification (\ref{V1.iden}), 
are equivalent to a vertex algebroid 
$(A,\Omega,\T,\ast,\{\;\},\{\;\}_\Omega)$.
The extended Lie algebroid $(A,\Omega,\T)$ consists of precisely those 
data that are independent of the choice of $s$, namely
\beq \label{VA.ELA}
\begin{array}{lll}
fg:=f_{(-1)}g &
f\alpha:=f_{(-1)}\alpha &
fX:=f_{(-1)}s(X)\;\t{mod}\;\Omega \vss \\
Xf:=s(X)_{(0)}f\ph{aa} &
L_X\alpha:=s(X)_{(0)}\alpha\ph{aa} &
[X,Y]:=s(X)_{(0)}s(Y)\;\t{mod}\;\Omega \vss \\
df:=Tf &
\alpha(X):=\alpha_{(1)}s(X) &
\end{array}
\eeq
\footnote{
For example, the definition of $Xf$ is indeed independent of $s$ because
$\alpha_{(0)}f=0$ for $f\in A$ and $\alpha\in\Omega$.
}
for $f,g\in A$, $\alpha\in\Omega$ and $X,Y\in\T$;
on the other hand
\begin{align} 
f\ast X\;\;\; &:= f_{(-1)}s(X)-s(fX) \nonumber \\
\{X,Y\}\;\; &:= s(X)_{(1)}s(Y) \label{VA.VAoid} \\
\{X,Y\}_\Omega &:= s(X)_{(0)}s(Y)-s([X,Y]) \nonumber
\end{align}
for $f\in A$ and $X,Y\in\T$.
\end{subsec}

\begin{subsec}
{\bf The induced morphism of vertex algebroids.}
Consider a homomorphism of vertex algebras $\Phi:V\rightarrow V'$.
Let $(A,\Omega,\T,\ast,\{\;\},\{\;\}_\Omega)$,\,
$(A',\Omega',\T',\ast',\{\;\}',\{\;\}'_\Omega)$
be the vertex algebroids associated to $V$, $V'$ and some splittings 
$s:\T\rightarrow V_1$, $s':\T'\rightarrow V'_1$.
The part of data of $\Phi$ involving only the two lowest weights, 
when rephrased in terms of identifications like (\ref{V1.iden}), 
are equivalent to a morphism $(\varphi,\Delta)$ between the two vertex
algebroids.
It consists of the obvious map of triples 
$\varphi:(A,\Omega,\T)\rightarrow(A',\Omega',\T')$ induced by $\Phi$, and 
a map $\Delta:\T\rightarrow\Omega'$ given by
\beqq
\Delta(X)=\Phi s(X)-s'(\varphi X),\quad X\in\T.
\eeqq
\end{subsec}

\begin{subsec} \label{VAoid.VA}
{\bf The vertex algebra freely generated by a vertex algebroid.}
Let $(A,\Omega,\T,\ast,\{\;\},\{\;\}_\Omega)$ be a vertex algebroid.
Throughout this discussion, we always have $f,g\in A$, 
$\alpha,\beta\in\Omega$, $X,Y\in\T$.
Define an associative $\bb{C}$-algebra $\mc{W}$ with generators of the form 
$f_n,\,\alpha_n,\,X_n,\,n\in\bb{Z}$ and the following relations
\beq \label{freeVA.comm}
\begin{array}{l}
(cf)_n=cf_n \qquad 
(c\alpha)_n=c\alpha_n \qquad 
(cX)_n=cX_n \vs \\
\mathbf 1_n=\delta_{n,0} \qquad\quad
(df)_n=-nf_n \qquad 
[f_n,g_m]=[f_n,\alpha_m]=[\alpha_n,\beta_m]=0 \vs \\ 
\lbrack X_n,f_m]=(Xf)_{n+m}\qquad\quad
[X_n,\alpha_m]=(L_X\alpha)_{n+m}+n\alpha(X)_{n+m} \vss \\
\lbrack X_n,Y_m]=[X,Y]_{n+m}+(\{X,Y\}_\Omega)_{n+m}+n(\{X,Y\})_{n+m}
\end{array} 
\eeq
where $c\in\bb{C}$, $n,m\in\bb{Z}$.
The subalgebra $\mc{W}_+\subset\mc{W}$ generated by $f_n$, $n>0$ and 
$\alpha_n,X_n$, $n\geq 0$ admits a trivial action on $\bb{C}$.
Let $\wt{V}:=\mc{W}\otimes_{\mc{W}_+}\bb{C}$ be the induced 
$\mc{W}$-module and $V:=\wt{V}/\sim$ the quotient module obtained by 
imposing the following relations for $v\in\wt{V}$:
\beq \label{freeVA.NOP}
\begin{array}{lcl}
(fg)_n v & \sim & \sum_{k\in\bb{Z}} f_k g_{n-k}v \vss \\
(f\alpha)_n v & \sim & \sum_{k\in\bb{Z}} f_k\alpha_{n-k}v \vss \\
(fX)_n v & \sim & \sum_{k\leq 0}f_k X_{n-k}v+\sum_{k>0}X_{n-k}f_k v
  -(f\ast X)_n v 
\end{array}
\eeq
Notice that the summations are always finite.
It is a consequence of the axioms of a vertex algebroid that 
(\ref{freeVA.comm})--(\ref{freeVA.NOP}) are consistent.
\footnote{
For example, $[X_n,(fY)_m]$ can be computed by either taking the commutator 
first or expanding $(fY)_m$ first.
The resulting identity is already implied by the vertex algebroid axioms
and does not lead to a new relation.
}
Define a vertex algebra structure on $V$ as follows.
The vacuum $\mathbf 1\in V$ is given by the coset of 
$1\otimes 1\in\wt{V}$.
The infinitesimal translation $T$ and weight operator $L_0$ are 
determined by the requirements
\beqq
\begin{array}{llll}
T\mathbf 1=0\ph{aaa} &
[T,f_n]=(1-n)f_{n-1}\ph{aa} &
[T,\alpha_n]=-n\alpha_{n-1}\ph{aa} &
[T,X_n]=-nX_{n-1} \vss \\
L_0\mathbf 1=0 &
[L_0,f_n]=-nf_n &
[L_0,\alpha_n]=-n\alpha_n &
[L_0,X_n]=-nX_n
\end{array}
\eeqq
which are consistent with (\ref{freeVA.comm})--(\ref{freeVA.NOP});
notice that actions of $f_n,\alpha_n,X_n$ change weights by $-n$.
Identify $f,\,\alpha,\,X$ with 
$f_0\mathbf 1,\,\alpha_{-1}\mathbf 1,\,X_{-1}\mathbf 1$ and associate
to them the following fields
\beqq
\begin{array}{c}
\sum_n f_n z^{-n}\,,\qquad
\sum_n \alpha_n z^{-n-1}\,,\qquad
\sum_n X_n z^{-n-1}
\end{array}
\eeqq
which are mutually local by (\ref{freeVA.comm});
notice that $f_{(n)}=f_{n+1}$, $\alpha_{(n)}=\alpha_n$, $X_{(n)}=X_n$.
Now apply the Strong Reconstruction Theorem \cite{FB-Z}.

Suppose $(A,\Omega,\T,\ast,\{\;\},\{\;\}_\Omega)$ is the vertex algebroid 
associated with a vertex algebra $V'$ and a splitting 
$s:\T\rightarrow V'_1$.
There is a canonical homomorphism of vertex algebras 
$\Phi:V\rightarrow V'$, determined by $\Phi f=f$, $\Phi\alpha=\alpha$ 
and $\Phi X=s(X)$.
If $\Phi$ is an isomorphism, $V'$ is said to be \emph{freely generated by
a vertex algebroid}. 
\end{subsec}

\begin{subsec} \label{VAoid.VA.mor}
{\bf The induced homomorphism of vertex algebras.}
A morphism of vertex algebroids
\beqq
(\varphi,\Delta):(A,\Omega,\T,\ast,\{\;\},\{\;\}_\Omega)\rightarrow
  (A',\Omega',\T',\ast',\{\;\}',\{\;\}'_\Omega)
\eeqq
induces a homomorphism $\Phi:V\rightarrow V'$ between the freely generated 
vertex algebras by the equations
\beqq
\begin{array}{lll}
\Phi f=\varphi f &
\Phi\alpha=\varphi\alpha & 
\Phi X=\varphi X+\Delta(X) \vss \\
\Phi\circ f_n=(\Phi f)_n\circ\Phi\qquad &
\Phi\circ\alpha_n=(\Phi\alpha)_n\circ\Phi\qquad &
\Phi\circ X_n=(\Phi X)_n\circ\Phi
\end{array}
\eeqq
for $f\in A,\,\alpha\in\Omega,\,X\in\T,n\in\bb{Z}$.
Indeed, these equations are consistent with 
(\ref{freeVA.comm})--(\ref{freeVA.NOP}).
\end{subsec}

\begin{subsec} \label{VAoid.VA.filtration}
{\bf More details on the constructions in \S\ref{VAoid.VA} and 
\S\ref{VAoid.VA.mor}.}
Given a possibly empty sequence of negative integers 
$\n=\{n_1\leq\cdots\leq n_s<0\}$, we write
\beqq
|\n|=n_1+\cdots+n_s\quad (0\t{ if }\n=\{\}),\qquad
\n(i)=\t{number of times }i\t{ appears in }\n
\eeqq
and regard $\n$ as a partition of $|\n|$.
For $k\geq 0$, let $I_k$ be the set of pairs $(\n,\m)$ of such sequences 
with $-|\n|-|\m|=k$.
Define a partial ordering on $I_k$ such that $(\n,\m)\prec(\n',\m')$ 
if and only if
\beqq
\begin{array}{lll}
-|\n|<-|\n'| &
\quad\t{or}\quad & 
|\n|=|\n'|\t{ and }\n'\t{ is a proper subpartition of }\n \vss \\
& \quad\t{or} & 
\n=\n'\t{ and }\m\t{ is a proper subpartition of }\m'
\end{array}
\eeqq
For example, $(\{\},\{-2,-2,-1\})\prec(\{\},\{-3,-2\})\prec(\{-4\},\{-1\})
\prec(\{-3,-1\},\{-1\})$ in $I_5$.

Consider the vertex algebra $V$ constructed in \S\ref{VAoid.VA}.
Associate to each $\n=\{n_1\leq\cdots\leq n_s<0\}$ and $s$-tuples 
$\bs{\alpha}=(\alpha_1,\ldots,\alpha_s)\in\Omega^s$, 
$\bs{X}=(X_1,\ldots,X_s)\in\T^s$ the following operators on $V$
\beqq
\bs{\alpha}_\n:=\alpha_{1,n_1}\cdots\alpha_{s,n_s},\qquad
\bs{X}_\n:=X_{1,n_1}\cdots X_{s,n_s}\qquad
(\t{both }1\t{ if }\n=\{\})
\eeqq
For $k>0$, we have 
$V_k=\t{span}\,\{\bs{X}_\n\bs{\alpha}_\m\bs{1}\,|\,(\n,\m)\in I_k\}$;
for $(\n,\m)\in I_k$, define the subspaces
\begin{align*}
\mc{F}_{\preceq(\n,\m)}&:=\t{span}\,
\{\bs{X}_{\n'}\bs{\alpha}_{\m'}\bs{1}\,|\,(\n',\m')\in I_k,\,
  (\n',\m')\preceq(\n,\m)\}\subset V_k \\
\mc{F}_{\prec(\n,\m)}&:=\t{span}\,
\{\bs{X}_{\n'}\bs{\alpha}_{\m'}\bs{1}\,|\,(\n',\m')\in I_k,\,
  (\n',\m')\prec(\n,\m)\}\subset V_k
\end{align*}
The bilinear operation $A\times V_k\rightarrow V_k$ given by 
$(f,v)\mapsto f_0 v$ does not make $V_k$ an $A$-module, but it preserves 
$\mc{F}_{\preceq(\n,\m)}$, $\mc{F}_{\prec(\n,\m)}$ and induces 
an $A$-module structure on their quotient.
In fact,
\beqq
\mc{F}_{\preceq(\n,\m)}/\mc{F}_{\prec(\n,\m)}\cong
\left(\bigotimes_{i=-1}^{-\infty}\t{Sym}_A^{\n(i)}\T\right)\otimes
\left(\bigotimes_{j=-1}^{-\infty}\t{Sym}_A^{\m(j)}\Omega\right)
\eeqq
as $A$-modules.
This allows us to compute the ``associated graded space''
\footnote{
More precisely, the coefficient of $q^k$ is the associated graded space 
of a certain filtration on $V_k$.
}
\beqq
A\oplus\bigoplus_{k=1}^\infty\left(q^k\bigoplus_{(\n,\m)\in I_k}
  \mc{F}_{\preceq(\n,\m)}/\mc{F}_{\prec(\n,\m)}\right)
\cong
\bigotimes_{\ell=1}^\infty\t{Sym}_{q^\ell}(\T\oplus\Omega) 
\eeqq
where $q$ is a formal variable and 
$\t{Sym}_t(\cdot)=\sum_{t=0}^\infty t^n\t{Sym}_A^n(\cdot)$.
The subspaces $\mc{F}_{\preceq(\n,\m)}$, $\mc{F}_{\prec(\n,\m)}$ and 
the isomorphisms stated here are natural, i.e.~respected by 
the homomorphism $\Phi$ constructed in \S\ref{VAoid.VA.mor}.
\end{subsec}

The omitted proofs of the following lemmas are straightforward 
(though somewhat tedious).

\begin{lemma} \label{lemma.newVAoid}
Given the following data: \vs \\
\indent $\cdot$\;
a vertex algebroid $(A,\Omega,\T,\ast,\{\;\},\{\;\}_\Omega)$ \vss \\
\indent $\cdot$\;
an isomorphism of extended Lie algebroids 
$\varphi:(A,\Omega,\T)\rightarrow(A',\Omega',\T')$ \vss \\
\indent $\cdot$\;
a $\bb{C}$-linear map $\Delta:\T\rightarrow\Omega'$ \vss \\
if we define maps
\beqq
\ast':A'\times\T'\rightarrow\Omega',\qquad
\{\;\}':\T'\times\T'\rightarrow A',\qquad
\{\;\}'_\Omega:\T'\times\T'\rightarrow\Omega'
\eeqq
by the equations in Definition \ref{VAoid.mor}, then 
$(A',\Omega',\T',\ast',\{\;\}',\{\;\}'_\Omega)$ is a vertex algebroid and
$(\varphi,\Delta)$ is by construction an isomorphism between the two 
vertex algebroids. \qedsymbol
\end{lemma}

\begin{lemma} \label{lemma.checkVAoidmor}
Given the following data: \vs \\
\indent $\cdot$\;
two vertex algebroids $(A,\Omega,\T,\ast,\{\;\},\{\;\}_\Omega)$ and
$(A',\Omega',\T',\ast',\{\;\}',\{\;\}'_\Omega)$ \vss \\
\indent $\cdot$\;
a morphism of extended Lie algebroids 
$\varphi:(A,\Omega,\T)\rightarrow(A',\Omega',\T')$ \vss \\
\indent $\cdot$\;
a $\bb{C}$-linear map $\Delta:\T\rightarrow\Omega'$ \vss \\
\indent $\cdot$\;
a subset $S\subset\T$ that is closed under $[\;]$ and spans $\T$ as 
an $A$-module \vss \\
if $(\varphi,\Delta)$ satisfies the equations in Definition \ref{VAoid.mor} 
for $(f,X,Y)\in A\times S^2$, then it also does for 
$(f,X,Y)\in A\times\T^2$ and hence is a morphism between the two given 
vertex algebroids. 
\qedsymbol
\end{lemma}

\begin{subsec}
{\bf Super version.}
There is no difficulty in generalizing the discussions in this appendix
to define \emph{extended Lie superalgebroids}, 
\emph{vertex superalgebroids}, and relate them to vertex superalgebras.
\end{subsec}

\begin{example} \label{VAoid.Lie}
{\bf The vertex algebroids associated to a Lie algebra.}
Consider a Lie algebra $\mf{g}$ over $\bb{C}$ and a vertex algebroid
of the form $(\bb{C},0,\mf{g},0,\lambda,0)$ with $\mf{g}$ acting
trivially on $\bb{C}$.
The second, fourth and last components are trivial by necessity.
The conditions on $\lambda:\mf{g}\times\mf{g}\rightarrow\bb{C}$ are
\beqq
\lambda(X,Y)=\lambda(Y,X),\qquad
\lambda([X,Y],Z)+\lambda(Y,[X,Z])=0
\eeqq
i.e.~it is a symmetric invariant bilinear form on $\mf{g}$.
Let
\beqq
V_\lambda(\mf{g})
  =\t{ the vertex algebra freely generated by }
  (\bb{C},0,\mf{g},0,\lambda,0).
\eeqq
In the case $\mf{g}$ is simple, finite-dimensional and $\lambda$
equals $k$ times the normalized Killing form, this is the vertex
algebra defined on the level-$k$ vacuum representation of
the affine Kac-Moody algebra $\hat{\mf{g}}$. \cite{FB-Z}
\end{example}

\begin{example} \label{algCDO}
{\bf Polynomial CDOs.}
Given nonnegative integers $p$ and $q$, let $\mc{W}$ be the associative 
$\bb{C}$-superalgebra generated by elements of the form
\beqq
b^i_n,\, a_{i,n},\qquad n\in\bb{Z},\; i=1,\ldots,p+q,\qquad
|b^i_n|=|a_{i,n}|=\begin{cases}
  \bar{0}, & i=1,\ldots,p \\
  \bar{1}, & i=p+1,\ldots,p+q
\end{cases}
\eeqq
($|\cdot|=$ parity) satisfying the following relations
\beqq
[b^i_n,b^j_m]=0=[a_{i,n},a_{j,m}],\qquad
[a_{i,n},b^j_m]=\delta^j_i\delta_{n,-m}
\eeqq
where $[\;]$ is the supercommutator.
The subalgebra $\mc{W}_+\subset\mc{W}$ generated by $b^i_n$, $n>0$ and 
$a_{i,n}$, $n\geq 0$ is supercommutative and admits a (purely even) 
trivial representation $\bb{C}$.
The induced $\mc{W}$-module
\beqq
\CDO(\bb{A}^{p|q}):=\mc{W}\otimes_{\mc{W}_+}\bb{C}
\eeqq
has the structure of a vertex superalgebra.
The vacuum is given by $\mathbf 1=1\otimes 1$.
The infinitesimal translation $T$ and weight operator $L_0$ are 
determined by
\beqq
\begin{array}{lll}
T\mathbf 1=0 &
[T,b^i_n]=(1-n)b^i_{n-1} \ph{aa} &
[T,a_{i,n}]=-na_{i,n-1} \vss \\
L_0\mathbf 1=0 \ph{aa} &
[L_0,b^i_n]=-nb^i_n &
[L_0,a_{i,n}]=-na_{i,n}
\end{array}
\eeqq
The vertex operators of $b^i_0\mathbf 1$ and $a_{i,-1}\mathbf 1$ are given 
respectively by the fields
\beqq
\begin{array}{l}
\sum_n b^i_n z^{-n},\qquad \sum_n a_{i,n}z^{-n-1}
\end{array}
\eeqq
while the other vertex operators follow from the Reconstruction Theorem 
\cite{FB-Z}.
\footnote{
This vertex superalgebra is the tensor product of $p$ copies of 
the $\beta\gamma$-system and $q$ copies of the $bc$-system.
}

The vertex superalgebra $\CDO(\bb{A}^{p|q})$ is freely generated by 
the associated vertex superalgebroid.
To describe the latter, consider the algebraic supermanifold
\beqq
\begin{array}{c}
\bb{A}^{p|q}:=\t{Spec}\left(\bb{C}[b^1,\cdots,b^d]\otimes
  \bigwedge(b^{p+1},\cdots,b^{p+q})\right)
\end{array}
\eeqq
and identify its functions, $1$-forms and vector fields with the following
subquotients of $\CDO(\bb{A}^{p|q})$: \vss \\
\indent $\cdot$\;
$\mc{O}(\bb{A}^{p|q})=\CDO(\bb{A}^{p|q})_0$ via $b^i=b^i_0\bs{1}$, 
$b^i b^j=b^i_0 b^j_0\bs{1}$, etc. \vss \\
\indent $\cdot$\;
$\Omega^1(\bb{A}^{p|q})\subset\CDO(\bb{A}^{p|q})_1$ via 
$db^i=b^i_{-1}\bs{1}$ \vss \\
\indent $\cdot$\;
$\T(\bb{A}^{p|q})=\CDO(\bb{A}^{p|q})_1/\Omega^1(\bb{A}^{p|q})$ via
$\d_i=\d/\d b^i=$ coset of $a_{i,-1}\bs{1}$ \vs \\
\footnote{
From another point of view, making these identifications dictates our 
(sign) conventions for calculus on $\bb{A}^{p|q}$ (or $\bb{R}^{p|q}$).
For example, it follows from $\alpha(X):=\alpha_{(1)}s(X)$ in (\ref{VA.ELA})
that $db^i(\d_j)=\e_j\delta^i_j$.
}
Then ``the'' vertex superalgebroid associated to $\CDO(\bb{A}^{p|q})$ is 
of the form 
\beqq
\big(\mc{O}(\bb{A}^{p|q}),\Omega^1(\bb{A}^{p|q}),\T(\bb{A}^{p|q}),
  \ast^c,\{\;\}^c,\{\;\}^c_\Omega\big).
\eeqq
The extended Lie superalgebroid structure consists of 
the usual differential on functions, 
Lie bracket on vector fields, 
Lie derivations by vector fields on functions and $1$-forms, and 
pairing between $1$-forms and vector fields.
Let $\e_i:=(-1)^{|b^i|}$.
If we use the splitting
\beqq
s:\T(\bb{A}^{p|q})\rightarrow\CDO(\bb{A}^{p|q})_1,\qquad
X=X^i\d_i\;\mapsto\;\e_i^{1+|X|}a_{i,-1}X^i
\eeqq
the rest of the vertex superalgebroid structure, as given by 
(\ref{VA.VAoid}), reads
\begin{align}
f\ast^c X\;\; &=  
  -(\e_i\e_j)^{1+|f|+|X|}(\d_j\d_i f)X^i db^j 
\nonumber \\
\{X,Y\}^c\; &= 
  -\e_j^{1+|X|+|Y|}(\d_j X^i)(\d_i Y^j)  
\label{CDO.VAoid456} \\
\{X,Y\}^c_\Omega &= 
  -(\e_j\e_k)^{1+|X|+|Y|}(\d_k\d_j X^i)(\d_i Y^j)db^k 
\nonumber
\end{align}
The superscript $c$ refers to the dependence on
coordinates.
\end{example}

\begin{example} \label{smoothCDO}
{\bf Local smooth CDOs.}
Let $b^1,\ldots,b^p$ and $b^{p+1},\ldots,b^{p+q}$ be respectively the real,
even and odd coordinates of $\bb{R}^{p|q}$, regarded as a smooth
cs-manifold, namely
\beqq
\textstyle
C^\infty(\bb{R}^{p|q})=C^\infty(\bb{R}^p)\otimes
  \bigwedge(b^{p+1},\ldots,b^{p+q})\otimes\bb{C}.
\eeqq
Let $\mb{W}$ be the restriction of $\bb{R}^{p|q}$ to an open set in 
$\bb{R}^p$.
Motivated by Example \ref{algCDO}, we define a vertex superalgebra
$\CDO(\mb{W})$ as follows.
The functions, $1$-forms and vector fields on $\mb{W}$ form an extended 
Lie superalgebroid as in Example \ref{algCDO}, and formulae 
(\ref{CDO.VAoid456}) again yield a vertex superalgebroid
\beqq
\big(C^\infty(\mb{W}),\Omega^1(\mb{W}),\T(\mb{W}),
  \ast^c,\{\;\}^c,\{\;\}^c_\Omega\big);
\eeqq
then take the freely generated vertex superalgebra. 
\end{example}

\begin{subsec} \label{MC.WZ}
{\bf The Wess-Zumino form of a diffeomorphism.}
Suppose $\bvarphi:\mb{W}\rightarrow\mb{W}'$ is a diffeomorphism
between restrictions of $\bb{R}^{p|q}$ (as a cs-manifold) to open sets.
The following notations will be used:
\beqq
|\cdot|=\t{parity},\quad
\e_i=(-1)^{|b^i|},\quad
\e_{ij}=(-1)^{|b^i||b^j|},\quad
i,j=1,\ldots,p+q
\eeqq
Let $g_\bvarphi:\mb{W}\rightarrow GL(p|q)$ be the map of
cs-manifolds whose components $(g_\bvarphi)^i_{\ph{i}j}$ are
given by
\beqq
\bvarphi^*db^i=(g_\bvarphi)^i_{\ph{i}j}db^j
\quad\Leftrightarrow\quad
(g_\bvarphi)^i_{\ph{i}j}=\e_j\e_{ij}\d_j\bvarphi^i
\eeqq
where $\bvarphi^i=\bvarphi^*b^i$.
\footnote{
In this notation, the chain rule reads 
$g_{\bvarphi'\bvarphi}=(\bvarphi^*g_{\bvarphi'})\cdot g_\bvarphi$.
}
Define the following differential forms
\beqq
\theta_\bvarphi:=g_\bvarphi^{-1}\cdot dg_\bvarphi
  \in\Omega^1(\mb{W})\otimes\mf{gl}(p|q), \qquad
WZ_\bvarphi:=\frac{1}{3}\Str(\theta_\bvarphi\wedge\theta_\bvarphi\wedge
  \theta_\bvarphi)\in\Omega^3(\mb{W}).
\eeqq
It follows from $d\theta_\bvarphi=-\theta_\bvarphi\wedge\theta_\bvarphi$ 
that $WZ_\bvarphi$ is closed.
\end{subsec}

\begin{theorem} \label{CDO.iso}
Let $\mb{W},\mb{W}',\mb{W}''$ be restrictions of $\bb{R}^{p|q}$ 
(as a cs-manifold) to open sets in $\bb{R}^p$. 

(a) Suppose $\bvarphi:\mb{W}\rightarrow\mb{W}'$ is a diffeomorphism.
There is a one-to-one correspondence:
\beqq
\left\{\hspace{-0.05in}
\begin{array}{c}
  \t{isomorphisms of vertex superalgebras }
  \CDO(\mb{W}')\rightarrow\CDO(\mb{W})\\ 
  \t{whose weight-zero components are }
  \bvarphi^*:C^\infty(\mb{W}')\rightarrow C^\infty(\mb{W})
\end{array}\hspace{-0.05in}
\right\}
\stackrel{\sim}{\longleftrightarrow}
\left\{\hspace{-0.05in}
\begin{array}{c}
  \xi\in\Omega^2(\mb{W}),\,\t{even}\\
  \t{and }d\xi=WZ_\bvarphi
\end{array}\hspace{-0.05in}
\right\}
\eeqq
Given $\xi$ as above, the corresponding isomorphism, denoted by 
$\bvarphi^*_\xi$, is induced by an isomorphism between the associated vertex
superalgebroids
\beqq
(\bvarphi^*,\Delta_{\bvarphi,\xi}):
\big(C^\infty(\mb{W}'),\Omega^1(\mb{W}'),\T(\mb{W}'),
  \ast^c,\{\;\}^c,\{\;\}^c_\Omega\big)
\rightarrow
\big(C^\infty(\mb{W}),\Omega^1(\mb{W}),\T(\mb{W}),
  \ast^c,\{\;\}^c,\{\;\}^c_\Omega\big)
\eeqq
where $\Delta_{\bvarphi,\xi}:\T(\mb{W}')\rightarrow\Omega^1(\mb{W})$ 
is given by
\beqq
\Delta_{\bvarphi,\xi}(X)=
-\e_i\e_{ij}\e_j^{1+|X|}\d_j(\bvarphi^*X)^i
  (\theta_\bvarphi)^j_{\ph{j}i}
-\frac{1}{2}\iota_{\bvarphi^*X}\Str(\theta_\bvarphi\otimes\theta_\bvarphi)
-\frac{1}{2}\iota_{\bvarphi^*X}\xi
\eeqq
for homogeneous elements.

(b) Suppose $\bvarphi':\mb{W}'\rightarrow\mb{W}''$ is another diffeomorphism,
$\xi'\in\Omega^2(\mb{W}')$ is even, and $d\xi'=WZ_{\bvarphi'}$.
Then we have the composition
\beqq
\bvarphi^*_\xi\circ\bvarphi^{\prime*}_{\xi'}
=(\bvarphi'\bvarphi)^*_{\eta},\qquad
\eta=\xi+\bvarphi^*\xi'+\sigma_{\bvarphi',\bvarphi}
\eeqq
where $\sigma_{\bvarphi',\bvarphi}:=\Str(\theta_\bvarphi\wedge
g_\bvarphi^{-1}\cdot\bvarphi^*\theta_{\bvarphi'}\cdot g_\bvarphi)$.
\end{theorem}

{\it Remarks.}
(i) This is a reformulation of a result in \cite{GMS1} in the smooth case.
(ii) As a consistency check, it follows from 
$\theta_{\bvarphi'\bvarphi}=\theta_\bvarphi
+g_\bvarphi^{-1}\cdot\bvarphi^*\theta_{\bvarphi'}\cdot g_\bvarphi$ that 
$WZ_{\bvarphi'\bvarphi}=WZ_\bvarphi+\bvarphi^*WZ_{\bvarphi'}+
d\sigma_{\bvarphi',\bvarphi}$.

\begin{proof}[Proof of Theorem \ref{CDO.iso}]
(a) Any morphism between the extended Lie superalgebroids associated to 
$\mb{W}$ and $\mb{W}'$ is induced by a map of cs-manifolds, which is 
$\bvarphi$ in this case.
Consider a morphism of vertex superalgebroids of the form
\beqq
(\bvarphi^*,\Delta):
\big(C^\infty(\mb{W}'),\Omega^1(\mb{W}'),\T(\mb{W}'),
  \ast^c,\{\;\}^c,\{\;\}^c_\Omega\big)\rightarrow
\big(C^\infty(\mb{W}),\Omega^1(\mb{W}),\T(\mb{W}),
  \ast^c,\{\;\}^c,\{\;\}^c_\Omega\big)
\eeqq
and write $\Delta=\Delta_0\circ\bvarphi^*|_{\T(\mb{W}')}$ in terms of 
an even map $\Delta_0:\T(\mb{W})\rightarrow\Omega^1(\mb{W})$.
Applying Lemma \ref{lemma.checkVAoidmor} with 
$S=\{\bvarphi_*\d_i\}_{i=1}^{p+q}$ and using (\ref{CDO.VAoid456}), we
obtain a complete set of equations for $\Delta_0$, namely
\beqq
\begin{array}{ll}
\Delta_0(f\d_i)-f\Delta_0(\d_i)
&=-\e_i\e_j^{1+|f|}(\d_j f)(\theta_\bvarphi)^j_{\ph{j}i} \vss \\
\Delta_0(\d_i)(\d_j)+\e_{ij}\Delta_0(\d_j)(\d_i)
&=-\Str(\theta_\bvarphi\otimes\theta_\bvarphi)(\d_i\otimes\d_j) \vss \\
\d_i\big(\Delta_0(\d_j)(\d_k)\big)
-\e_{ij}\d_j\big(\Delta_0(\d_i)(\d_k)\big)
+\e_{ik}\e_{jk}\d_k\big(\Delta_0(\d_i)(\d_j)\big) \hspace{-1.6in} & \vss \\
&=-\e_r\e_{jk}\e_{kr}\e_{ks}\,
  \d_i\big((\theta_\bvarphi)^r_{\ph{i}s}(\d_k)\big)\cdot
  (\theta_\bvarphi)^s_{\ph{i}r}(\d_j)
\end{array}
\eeqq
The first equation implies that for any $X\in\T(\mb{W})$
\beqq
\Delta_0(X)=
-\e_i\e_{ij}\e_j^{1+|X|}(\d_j X^i)(\theta_\bvarphi)^j_{\ph{j}i}
+X^i\Delta_0(\d_i)\,;
\eeqq
the second equation allows us to write
\beqq
\Delta_0(\d_i)(\d_j)=
-\frac{1}{2}\Str(\theta_\bvarphi\otimes\theta_\bvarphi)(\d_i\otimes\d_j)
-\frac{1}{2}\xi_{ij},\qquad
\xi_{ji}=-\e_{ij}\xi_{ij}\,;
\eeqq
then it follows from
$d\theta_\bvarphi=-\theta_\bvarphi\wedge\theta_\bvarphi$ that the third
equation is equivalent to
\beqq
d\xi=WZ_\bvarphi
=\frac{1}{3}\Str(\theta_\bvarphi\wedge\theta_\bvarphi\wedge\theta_\bvarphi)
\eeqq
where $\xi$ is the even $2$-form with $\xi(\d_i,\d_j)=\xi_{ij}$.
Since $\CDO(\mb{W})$ and $\CDO(\mb{W}')$ are freely generated by vertex 
superalgebroids, an isomorphism between them is equivalent to 
an isomorphism between the associated vertex superalgebroids.
This completes the proof of (a).

(b) By part (a), the composition in question must be of the form 
$(\bvarphi'\bvarphi)^*_\eta$ for some $\eta\in\Omega^2(\mb{W})$.
At the level of vertex superalgebroids, we have
\begin{align*}
(\bvarphi^*\bvarphi^{\prime*},\Delta_{\bvarphi'\bvarphi,\eta})
&=(\bvarphi^*,\Delta_{\bvarphi,\xi})\circ
  (\bvarphi^{\prime*},\Delta_{\bvarphi',\xi'}) \\
\Leftrightarrow\qquad\qquad
\Delta_{\bvarphi'\bvarphi,\eta}
&=\bvarphi^*\Delta_{\bvarphi',\xi'}+\Delta_{\bvarphi,\xi}\bvarphi^{\prime*}.
\end{align*}
Evaluation at e.g.~$\bvarphi'_*\bvarphi_*\d_k$ then yields the desired
formula for $\eta$.
\end{proof}

\newpage

\setcounter{equation}{0}
\section{Affine Connections on CS-Manifolds}
\label{app.oddE.conn}

Consider a smooth manifold $M$ and a smooth $\bb{C}$-vector
bundle $E\rightarrow M$.
In this appendix, we construct an affine connection on
the smooth cs-manifold $\mb{M}=\Pi E$ and obtain a number of
formulae used in the computations of CDOs on $\mb{M}$.

\begin{subsec}
{\bf Functions on $\mb{M}$.}
Let $d=\dim M$ and $r=\t{rank}\,E$.
There is a canonical identification
\beq \label{cs.fcns}
C^\infty(\mb{M})\cong\Gamma(\wedge^*E^\vee).
\eeq
In particular, a set of local coordinates $(x^1,\cdots,x^d)$
on $M$ and a local frame $(\ve^1,\cdots,\ve^r)$ of $E^\vee$
together determine a set of local coordinates
$(x^1,\cdots,x^d,\ve^1,\cdots,\ve^r)$ on $\mb{M}$.
\end{subsec}

\begin{subsec} \label{sec.cs.vf}
{\bf Vector fields on $\mb{M}$.}
Choose a connection $\nabla^E$ on $E$ and use the same
notation for the induced connection on $\wedge^*E^\vee$.
Let $(\ve_1,\cdots,\ve_r)$ be the local frame of $E$
dual to $(\ve^1,\cdots,\ve^r)$.
Let $X,Y\in\T(M)$ and $\sigma,\tau\in\Gamma(E)$.
Under the identification (\ref{cs.fcns}), vector fields
on $\mb{M}$ correspond to derivations on sections of
$\wedge^*E^\vee$.
In particular, denote by
\beq \label{cs.vf}
\left\{\begin{array}{l}
  \mc{D}_X \\
  \mc{I}_\sigma \\
  J
\end{array}\right\}
\begin{array}{c}
  \t{the vector fields on }\mb{M} \\
  \t{corresponding to the }
\end{array}
\left\{\begin{array}{l}
  \t{covariant differentiation }\nabla^E_X \\
  \t{contraction with }\sigma \\
  \t{exterior degree}
\end{array}\right\}
\eeq 
The vector fields $\mc{D}_X$ and $\mc{I}_\sigma$
span $\T(\mb{M})$ over $C^\infty(\mb{M})$.
The super Lie brackets of (\ref{cs.vf}) are given by
\beq \label{cs.bracket}
[\mc{D}_X,\mc{D}_Y]
  =\mc{D}_{[X,Y]}
  -\ve^k\mc{I}_{R^E_{X,Y}\ve_k},\quad
[\mc{D}_X,\mc{I}_\sigma]
  =\mc{I}_{\nabla^E_X\sigma},\quad
[\mc{I}_\sigma,\mc{I}_\tau]=0=[J,\mc{D}_X],\quad
[J,\mc{I}_\sigma]=-\mc{I}_\sigma
\eeq
where $R^E$ is the curvature of $\nabla^E$.
\end{subsec}

\begin{subsec} \label{sec.cs.conn}
{\bf An affine connection on $\mb{M}$.}
Choose also a connection $\nabla^M$ on $TM$.
Let $X,Y,Z\in\T(M)$ and $\sigma,\tau\in\Gamma(E)$.
Define a connection $\nabla$ on $T\mb{M}$ by
\beq \label{cs.conn}
\nabla_{\mc{D}_X}\mc{D}_Y=\mc{D}_{\nabla^M_X Y},\quad
\nabla_{\mc{D}_X}\mc{I}_\sigma=\mc{I}_{\nabla^E_X\sigma},\quad
\nabla_{\mc{I}_\sigma}\mc{D}_X
  =\nabla_{\mc{I}_\sigma}\mc{I}_\tau=0
\eeq
and the Leibniz rule.
Using (\ref{cs.bracket}), we compute the curvature
of $\nabla$ as follows
\beq \label{cs.curv}
R_{\mc{D}_X,\mc{D}_Y}\mc{D}_Z=\mc{D}_{R^M_{X,Y}Z},\quad
R_{\mc{D}_X,\mc{D}_Y}\mc{I}_\sigma=\mc{I}_{R^E_{X,Y}\sigma},\quad
R_{\mc{D}_X,\mc{I}_\sigma}=R_{\mc{I}_\sigma,\mc{I}_\tau}=0
\eeq
where $R^M$ is the curvature of $\nabla^M$.
\end{subsec}

\begin{lemma} \label{lemma.ops.J}
(a) The operator $\tnabla J$ sends $\mc{D}_X$ to $0$,
and $\mc{I}_\sigma$ to itself.
(b) $\nabla(\tnabla J)=0$.
\end{lemma}

\begin{proof}
Recall that $\tnabla J:=\nabla_J-[J,\,\t{-}\,]$.
Using the fact that $J=\ve^k\mc{I}_{\ve_k}$,
(a) follows readily from (\ref{cs.bracket}) and
(\ref{cs.conn}).
Then (b) is clear.
\end{proof}

\begin{lemma} \label{lemma.str}
Regarding $\Omega^*(M)$ as a subalgebra of $\Omega^*(\mb{M})$,
we have

(a) $\Str R=\Tr R^M-\Tr R^E$,

(b) $\Str(R\wedge R)=\Tr(R^M\wedge R^M)-\Tr(R^E\wedge R^E)$,

(c) $\Str(R\cdot\tnabla J)=-\Tr R^E$.
\end{lemma}

\begin{proof}
All these statements follow easily from (\ref{cs.curv})
and Lemma \ref{lemma.ops.J}a.
\end{proof}

\begin{example} \label{dRcs}
{\bf The de Rham cs-manifold.}
In the case $E=TM\otimes\bb{C}$, (\ref{cs.fcns}) can be
rewritten as
\beqq
C^\infty(\mb{M})\cong\Omega^*(M).
\eeqq
Let $\ve^i=dx^i$ and $\ve_i=\d_i=\d/\d x^i$.
Besides (\ref{cs.vf}), consider also the odd vector field
$Q=\ve^i\d_i$ on $\mb{M}$ corresponding to the de Rham operator $d$.
Assume that $\nabla^M$ is torsion-free.
This implies the identity $d=dx^i\wedge\nabla^M_{\d_i}$,
or equivalently
\beq \label{Q.DQ}
Q=\ve^i\mc{D}_{\d_i}=\mc{D}_Q
\eeq
where the second equality should be understood as
the definition of a new notation.
Similar abuse of notation will appear below without
further comment.
The super Lie brackets with $Q$ are given by
\beq \label{Q.bracket}
[Q,\mc{D}_X]
  =\mc{D}_{\nabla^M_Q X}
  -\mc{I}_{R^M_{Q,X}Q},\quad
[Q,\mc{I}_X]
  =\mc{I}_{\nabla^M_Q X}+\mc{D}_X,\quad
[J,Q]=Q,\quad
[Q,Q]=0.
\eeq
Indeed, the first two equations follow from the following
calculations
\begin{align*}
[Q,\mc{D}_X]
&=[\ve^i\mc{D}_{\d_i},\mc{D}_X]
 =\ve^i[\mc{D}_{\d_i},\mc{D}_X]-(\mc{D}_X\ve^i)\mc{D}_{\d_i}
 =\ve^i\mc{D}_{[\d_i,X]}
  -\ve^i\ve^j\mc{I}_{R^M_{\d_i,X}\d_j}
  +\ve^i\mc{D}_{\nabla^M_X\d_i} \\
&\textstyle
=\ve^i\mc{D}_{\nabla^M_{\d_i}X}-\mc{I}_{R^M_{Q,X}Q}
=\mc{D}_{\nabla^M_Q X}-\mc{I}_{R^M_{Q,X}Q} \\
[Q,\mc{I}_X]
&=[\ve^i\mc{D}_{\d_i},\mc{I}_X]
 =\ve^i[\mc{D}_{\d_i},\mc{I}_X]+(\mc{I}_X\ve^i)\mc{D}_{\d_i}
 =\ve^i\mc{I}_{\nabla^M_{\d_i}X}+\mc{D}_X
 =\mc{I}_{\nabla^M_Q X}+\mc{D}_X
\end{align*}
where we have used (\ref{Q.DQ}), (\ref{cs.bracket}) and
the torsion-free condition.
By (\ref{cs.conn}) and (\ref{Q.DQ}), covariant
differentiation with respect to $Q$ is given by
\beq \label{Q.conn}
\nabla_Q\mc{D}_X=\mc{D}_{\nabla^M_Q X},\quad
\nabla_Q\mc{I}_X=\mc{I}_{\nabla^M_Q X}.
\eeq
\end{example}

\begin{lemma} \label{lemma.DtQ}
The operator $\tnabla Q$ and its covariant derivatives
are computed as follows:

(a) $\tnabla Q$ sends
$\mc{D}_X$ to $\mc{I}_{R^M_{Q,X}Q}$, and
$\mc{I}_X$ to $-\mc{D}_X$.

(b) $\nabla_{\mc{D}_X}(\tnabla Q)$ sends
$\mc{D}_Y$ to $\mc{I}_{(\nabla^M_X R^M)_{Q,Y}Q}$, and
$\mc{I}_Y$ to $0$.

(c) $\nabla_{\mc{I}_X}(\tnabla Q)$ sends
$\mc{D}_Y$ to $\mc{I}_{R^M_{X,Q}Y}$, and
$\mc{I}_Y$ to $0$.
\end{lemma}

\begin{proof}
Recall that $\tnabla Q:=\nabla_Q-[Q,\,\t{-}\,]$.
(a) follows from (\ref{Q.bracket}) and (\ref{Q.conn}).
For (b) and (c) we compute 
\begin{align*}
\Big(\nabla_{\mc{D}_X}(\tnabla Q)\Big)\mc{D}_Y
&\textstyle
=\nabla_{\mc{D}_X}\mc{I}_{R^M_{Q,Y}Q}
  -(\tnabla Q)\mc{D}_{\nabla^M_X Y}
=\mc{I}_{\nabla^M_X R^M_{Q,Y}Q}
  +(\mc{D}_X\ve^i\ve^j)\mc{I}_{R^M_{\d_i,Y}\d_j}
  -\mc{I}_{R^M_{Q,\nabla^M_X Y}Q} \\
&\textstyle
=\mc{I}_{\nabla^M_X R^M_{Q,Y}Q}
  -\ve^i\mc{I}_{R^M_{\nabla^M_X\d_i,Y}Q}
  -\ve^i\mc{I}_{R^M_{Q,Y}\nabla^M_X\d_i}
  -\mc{I}_{R^M_{Q,\nabla^M_X Y}Q}
=\mc{I}_{(\nabla^M_X R^M)_{Q,Y}Q} \\
\Big(\nabla_{\mc{D}_X}(\tnabla Q)\Big)\mc{I}_Y
&=-\nabla_{\mc{D}_X}\mc{D}_Y
  -(\tnabla Q)\mc{I}_{\nabla^M_X Y}
=-\mc{D}_{\nabla^M_X Y}+\mc{D}_{\nabla^M_X Y}=0 \\
\Big(\nabla_{\mc{I}_X}(\tnabla Q)\Big)\mc{D}_Y
&=\nabla_{\mc{I}_X}\mc{I}_{R^M_{Q,Y}Q}
=(\mc{I}_X\ve^i\ve^j)\mc{I}_{R^M_{\d_i,Y}\d_j}
=\mc{I}_{R^M_{X,Y}Q}-\mc{I}_{R^M_{Q,Y}X}
=\mc{I}_{R^M_{X,Q}Y} \\
\Big(\nabla_{\mc{I}_X}(\tnabla Q)\Big)\mc{I}_Y
&=-\nabla_{\mc{I}_X}\mc{D}_Y=0
\end{align*}
where we have used (\ref{cs.conn}) and the first Bianchi
identity.
\end{proof}

\begin{lemma} \label{lemma.str.Q}
The operators $\tnabla Q$, $R\cdot\tnabla Q$,
$R\cdot\tnabla Q\cdot\tnabla Q$ and
$\nabla(\tnabla Q)\wedge\nabla(\tnabla Q)$
all have supertrace zero.
It follows that the supertrace of
$\nabla(\tnabla Q)\cdot\tnabla Q$ is closed.
\end{lemma}

\begin{proof}
The first three operators have no supertrace
by Lemma \ref{lemma.DtQ}a and (\ref{cs.curv}).
For the third, notice that
\beqq
\textstyle
(\tnabla Q)^2\mc{D}_X=-\mc{D}_{R^M_{Q,X}Q},\quad
(\tnabla Q)^2\mc{I}_X=-\mc{I}_{R^M_{Q,X}Q}.
\eeqq
The fourth operator has no supertrace by
Lemma \ref{lemma.DtQ}b and c.
The remaining assertion follows from
\beqq
d\Str\Big(\nabla(\tnabla Q)\cdot\tnabla Q\Big)
=2\Str(R\cdot\tnabla Q\cdot\tnabla Q)
-\Str\Big(\nabla(\tnabla Q)\wedge\nabla(\tnabla Q)\Big).\qedhere
\eeqq
\end{proof}

\begin{example} \label{Dolbcs}
{\bf Dolbeault cs-manifolds.}
Now we change our notations as follows:~$M$ is
a complex manifold, $TM$ its holomorphic tangent bundle,
$E$ a holomorphic vector bundle over $M$, and
$\mb{M}=\Pi(\overline{TM}\oplus E)$ as a smooth cs-manifold.
There is a canonical identification
\beq \label{Dolb.fcns}
C^\infty(\mb{M})\cong\Omega^{0,*}(M;\wedge^*E^\vee).
\eeq
Given a set of local holomorphic coordinates
$(z^1,\cdots,z^d)$ on $M$ and a local holomorphic frame
$(\ve^1,\cdots,\ve^r)$ of $E^\vee$, there is an associated set
of local coordinates on $\mb{M}$, namely
\beqq
(\t{Re}z^1,\t{Im}z^1,\cdots,\t{Re}z^d,\t{Im}z^d,
  \bar\zeta^1,\cdots,\bar\zeta^d,
  \ve^1,\cdots,\ve^r)
\eeqq
where $\bar\zeta^i$ correspond to $d\bar z^i$ under
(\ref{Dolb.fcns}).
Let $\dbar_i=\d/\d\bar z^i$ and $(\ve_1,\cdots,\ve_k)$ be the dual local
frame of $E$.

Choose connections $\nabla^M$ on $TM$ and $\nabla^E$ on $E$
of type $(1,0)$;
denote by $\bar\nabla^M$ the induced connection on
$\overline{TM}$.
Let $X,Y,Z\in\T(M)$, $U,V\in\T^{0,1}(M)$ and $\sigma,\tau\in\Gamma(E)$.
Under the identification (\ref{Dolb.fcns}), vector fields
on $\mb{M}$ correspond to derivations of $(0,*)$-forms on $M$
valued in $\wedge^*E^\vee$.
In particular, denote by
\beq \label{Dolb.vf}
\left\{\begin{array}{l}
  \mc{D}_X \\
  \mc{I}_U, \mc{I}_\sigma \\
  J^r, J^\ell \\
  Q
\end{array}\right\}
\begin{array}{c}
  \t{the vector fields on }\mb{M} \\
  \t{corresponding to the }
\end{array}
\left\{\begin{array}{l}
  \t{covariant differentiation }
    \bar\nabla^M_X\otimes 1+1\otimes\nabla^E_X \\
  \t{contractions with }U, \sigma \\
  \t{exterior degrees in }
    \wedge^*\overline{TM}{}^\vee, \wedge^*E^\vee \\
  \t{Dolbeault operator }\dbar\otimes 1
\end{array}\right\}
\eeq
The vector fields $\mc{D}_X$, $\mc{I}_U$ and $\mc{I}_\sigma$
span $\T(\mb{M})$ over $C^\infty(\mb{M})$.
Adopt an abuse of notation similar to that
in Example \ref{dRcs}.
The super Lie brackets among the first three types of vector
fields in (\ref{Dolb.vf}) are
\beq
\nonumber
&[\mc{D}_X,\mc{D}_Y]
  =\mc{D}_{[X,Y]}
  -\mc{I}_{\bar R^M_{X,Y}Q}
  -\ve^k\mc{I}_{R^E_{X,Y}\ve_k},\quad
[\mc{D}_X,\mc{I}_U]
  =\mc{I}_{\bar\nabla^M_X U},\quad
[\mc{D}_X,\mc{I}_\sigma]
  =\mc{I}_{\nabla^E_X\sigma}& \\
\label{Dolb.bracket}
&[\mc{I}_U,\mc{I}_V]
=[\mc{I}_U,\mc{I}_\sigma]
=[\mc{I}_\sigma,\mc{I}_\tau]
=0
=[J^r,\mc{D}_X]
=[J^r,\mc{I}_\sigma]
=[J^\ell,\mc{D}_X]
=[J^\ell,\mc{I}_U]
=[J^r,J^\ell] & \\
\nonumber
&[J^r,\mc{I}_U]=-\mc{I}_U,\quad
[J^\ell,\mc{I}_\sigma]=-\mc{I}_\sigma
\eeq
Assume that $\nabla^M$ is torsion-free.
\footnote{ \label{torsionfree}
If $\nabla^M$ has a nontrivial torsion $T$, we can
replace it with a new connection $\nabla^{\prime M}$ defined by
$\nabla^{\prime M}_X=\nabla^M_X-\frac{1}{2}T_{X,\,\t{-}\,}$,
which is also of type $(1,0)$ and is torsion-free.
}
This implies the identity
$\dbar=d\bar z^i\wedge\bar\nabla^M_{\dbar_i}$,
or equivalently
\beq \label{DolbQ.DQ}
Q=\bar\zeta^i\mc{D}_{\dbar_i}=\mc{D}_Q.
\eeq
Then the various super Lie brackets with $Q$ are given by
\beq \label{DolbQ.bracket}
\begin{array}{c}
[Q,\mc{D}_X]
  =\mc{D}_{\nabla^M_Q X^{1,0}+\bar\nabla^M_Q X^{0,1}}
  -\mc{I}_{\bar R^M_{Q,X}Q}
  +\ve^k\mc{I}_{R^E_{Q,X}\ve_k},\quad
[Q,\mc{I}_U]
  =\mc{I}_{\bar\nabla^M_Q U}+\mc{D}_U \vss \\
\,[Q,\mc{I}_\sigma]=\mc{I}_{\nabla^E_Q\sigma},\quad
[J^r,Q]=Q,\quad
[J^\ell,Q]
=0=[Q,Q]
\end{array}
\eeq
Indeed, the first two follow from (\ref{DolbQ.DQ}) and
calculations similar to those below (\ref{Q.bracket}).

Define a connection $\nabla$ on $T\mb{M}$ as in
\S\ref{sec.cs.conn}.
More explicitly, we define
\beq \label{Dolb.conn}
\begin{array}{c}
\nabla_{\mc{D}_X}\mc{D}_Y
  =\mc{D}_{\nabla^M_X Y^{1,0}+\bar\nabla^M_X Y^{0,1}},\quad
\nabla_{\mc{D}_X}\mc{I}_U
  =\mc{I}_{\bar\nabla^M_X U},\quad
\nabla_{\mc{D}_X}\mc{I}_\sigma
  =\mc{I}_{\nabla^E_X\sigma} \vss\\
\nabla_{\mc{I}_U}\mc{D}_X
=\nabla_{\mc{I}_U}\mc{I}_V
=\nabla_{\mc{I}_U}\mc{I}_\sigma
=\nabla_{\mc{I}_\sigma}\mc{D}_X
=\nabla_{\mc{I}_\sigma}\mc{I}_U
=\nabla_{\mc{I}_\sigma}\mc{I}_\tau
=0
\end{array}
\eeq
By (\ref{DolbQ.DQ}), covariant differentiation with respect to
$Q$ is given by
\beq \label{DolbQ.conn}
\nabla_Q\mc{D}_X
  =\mc{D}_{\nabla^M_Q X^{1,0}+\bar\nabla^M_Q X^{0,1}},\quad
\nabla_Q\mc{I}_U
  =\mc{I}_{\bar\nabla^M_Q U},\quad
\nabla_Q\mc{I}_\sigma
  =\mc{I}_{\nabla^E_Q\sigma}.
\eeq
Using (\ref{Dolb.bracket}), we compute the curvature
of $\nabla$ as follows
\beq \label{Dolb.curv}
\begin{array}{c}
R_{\mc{D}_X,\mc{D}_Y}\mc{D}_Z
  =\mc{D}_{R^M_{X,Y}Z^{1,0}+\bar R^M_{X,Y}Z^{0,1}},\quad
R_{\mc{D}_X,\mc{D}_Y}\mc{I}_U
  =\mc{I}_{\bar R^M_{X,Y}U},\quad
R_{\mc{D}_X,\mc{D}_Y}\mc{I}_\sigma
  =\mc{I}_{R^E_{X,Y}\sigma} \vss \\
R_{\mc{D}_X,\mc{I}_U}
=R_{\mc{D}_X,\mc{I}_\sigma}
=R_{\mc{I}_U,\mc{I}_V}
=R_{\mc{I}_U,\mc{I}_\sigma}
=R_{\mc{I}_\sigma,\mc{I}_\tau}
=0
\end{array}
\eeq
\end{example}

The following statements and their proofs are similar to
Lemmas \ref{lemma.ops.J} and \ref{lemma.str}.

\begin{lemma} \label{lemma.Jr}
(a) The operator $\tnabla J^r$ sends
$\mc{D}_X$, $\mc{I}_\sigma$ to $0$, and
$\mc{I}_U$ to itself.
(b) $\nabla(\tnabla J^r)=0$.
\end{lemma}

\begin{proof}
Use (\ref{Dolb.bracket}), (\ref{Dolb.conn}) and
the fact that $J^r=\bar\zeta^i\mc{I}_{\dbar_i}$.
\end{proof}

\begin{lemma} \label{lemma.Jl}
(a) The operator $\tnabla J^\ell$ sends
$\mc{D}_X$, $\mc{I}_U$ to $0$, and
$\mc{I}_\sigma$ to itself.
(b) $\nabla(\tnabla J^\ell)=0$.
\end{lemma}

\begin{proof}
Use (\ref{Dolb.bracket}), (\ref{Dolb.conn}) and
the fact that $J^\ell=\ve^k\mc{I}_{\ve_k}$.
\end{proof}

\begin{lemma} \label{lemma.str.Dolb}
Regarding $\Omega^*(M)$ as a subalgebra of
$\Omega^*(\mb{M})$, we have

(a) $\Str R=\Tr R^M-\Tr R^E$,

(b) $\Str(R\wedge R)=\Tr(R^M\wedge R^M)-\Tr(R^E\wedge R^E)$,

(c) $\Str(R\cdot\tnabla J^r)=-\Tr\bar R^M$ and
$\Str(R\cdot\tnabla J^\ell)=-\Tr R^E$.
\end{lemma}

\begin{proof}
Use (\ref{Dolb.curv}) and the previous two lemmas.
\end{proof}

The following statements and their proofs are similar
to Lemmas \ref{lemma.DtQ} and \ref{lemma.str.Q}.

\begin{lemma}
The operator $\tnabla Q$ and its covariant derivatives
are computed as follows:

(a) $\tnabla Q$ sends
$\mc{D}_X$ to $\mc{I}_{\bar R^M_{Q,X}Q}-\ve^k\mc{I}_{R^E_{Q,X}\ve_k}$,
$\mc{I}_U$ to $-\mc{D}_U$, and
$\mc{I}_\sigma$ to $0$.

(b) $\nabla_{\mc{D}_X}(\tnabla Q)$ sends $\mc{D}_Y$ to
$\mc{I}_{(\bar\nabla^M_X\bar R^M)_{Q,Y}Q}
-\ve^k\mc{I}_{(\nabla^E_X R^E)_{Q,Y}\ve_k}$,
and $\mc{I}_U$, $\mc{I}_\sigma$ to $0$.

(c) $\nabla_{\mc{I}_U}(\tnabla Q)$ sends
$\mc{D}_X$ to $\mc{I}_{\bar R^M_{U,Q}X}+\ve^k\mc{I}_{R^E_{U,X}\ve_k}$,
and $\mc{I}_V$, $\mc{I}_\sigma$ to $0$.

(d) $\nabla_{\mc{I}_\sigma}(\tnabla Q)$ sends
$\mc{D}_X$ to $-\mc{I}_{R^E_{Q,X}\sigma}$,
and $\mc{I}_U$, $\mc{I}_\sigma$ to $0$.
\end{lemma}

\begin{proof}
Use (\ref{DolbQ.bracket}), (\ref{Dolb.conn}),
(\ref{DolbQ.conn}) and the first Bianchi identity.
\end{proof}

\begin{lemma} \label{lemma.str.DolbQ}
The operators $\tnabla Q$, $R\cdot\tnabla Q$,
$R\cdot\tnabla Q\cdot\tnabla Q$ and
$\nabla(\tnabla Q)\wedge\nabla(\tnabla Q)$
all have supertrace zero.
It follows that the supertrace of
$\nabla(\tnabla Q)\cdot\tnabla Q$ is closed.
\end{lemma}

\begin{proof}
Use (\ref{Dolb.curv}) and the previous lemma.
For the third operator, also notice that
$\bar R^M_{Q,X}Q=\frac{1}{2}\bar R^M_{Q,Q}X^{0,1}$
by the first Bianchi identity,
and $R^E_{Q,U}=0$ by our assumption on $\nabla^E$.
\end{proof}

\newpage

{\footnotesize



}

\end{document}